\pdfoutput=1
\documentclass[numbers]{article}

\usepackage{arxiv}
\usepackage[utf8]{inputenc} % allow utf-8 input
\usepackage[T1]{fontenc}    % use 8-bit T1 fonts
\usepackage{hyperref}       % hyperlinks
\usepackage{url}            % simple URL typesetting
\usepackage{amsfonts}       % blackboard math symbols
\usepackage{algorithm,algorithmic}
\usepackage{amsmath}
\usepackage{amsthm}
\usepackage{listings}
\usepackage{xcolor}
\usepackage{a4wide}
\usepackage{latexsym}
\usepackage{graphicx}
\usepackage{subcaption}
\usepackage{varioref}
\usepackage{color}
\usepackage{lscape}
\usepackage{verbatim}
\usepackage{indentfirst}
\usepackage[section]{placeins}
\usepackage[abs]{overpic}
\usepackage{multirow}
\usepackage{caption}

\definecolor{gray}{rgb}{0.4,0.4,0.4}
\definecolor{darkblue}{rgb}{0.0,0.0,0.6}
\definecolor{cyan}{rgb}{0.0,0.6,0.6}
\definecolor{light-gray}{gray}{0.85}
\definecolor{mygreen}{rgb}{0,0.4,0}
\definecolor{mygray}{rgb}{0.5,0.5,0.5}
\definecolor{mylilas}{RGB}{170,55,241}

\newcommand{\calA}{{\cal A}} \newcommand{\calB}{{\cal B}}
\newcommand{\calC}{{\cal C}} 
 \newcommand{\calF}{{\cal F}}

 \newcommand{\calL}{{\cal L}}
\newcommand{\calM}{{\cal M}} \newcommand{\calN}{{\cal N}}
 \newcommand{\calP}{{\cal P}}
 \newcommand{\calR}{{\cal R}}
\newcommand{\calS}{{\cal S}} \newcommand{\calT}{{\cal T}}
\newcommand{\calU}{{\cal U}} 
 
\newcommand{\calY}{{\cal Y}} \newcommand{\calZ}{{\cal Z}}

\newcommand{\eqdef}{\stackrel{\rm def}{=}}
\newcommand{\sfrac}[2]{{\scriptstyle \frac{#1}{#2}}}
\newcommand{\half}{\sfrac{1}{2}}
\newcommand{\ip}[2]{\langle #1, #2 \rangle}
\newcommand{\tim}[1]{\;\; \mbox{#1} \;\;}

\makeatletter
\renewcommand{\paragraph}{\@startsection{paragraph}{4}{0ex}%
    {-3.25ex plus -1ex minus -0.2ex}%
    {1.5ex plus 0.2ex}%
    {\normalfont\normalsize\bfseries}}
\makeatother
 
\stepcounter{secnumdepth}
\stepcounter{tocdepth}

\newcounter{algor}[section]
\renewcommand{\thealgor}{\thesection.\arabic{algor}}

\newcommand{\algor}[3]{\refstepcounter{algor}
{\vspace{0,7cm}
\hrule 
{\raggedleft {\bf Algorithm \thealgor : #2}\label{#1} }
\vspace{0,1cm}
\hrule
{#3} 
\hrule
\vspace{0,7cm} }
}

\theoremstyle{definition}
\newtheorem{definition}{Definition}[section]
\title{DEFT-FUNNEL: an open-source global optimization solver for constrained grey-box and \\ black-box problems}

\author{
  Phillipe R.~Sampaio \\
  Veolia Research and Innovation\\
  Maisons-Laffitte, France \\
  \texttt{sampaio.phillipe@gmail.com}
}

\begin{document}
\maketitle

\begin{abstract}
The fast-growing need for grey-box and black-box optimization methods for constrained global optimization problems in fields such as medicine, chemistry, engineering and artificial intelligence, has contributed for the design of new efficient algorithms for finding the best possible solution. In this work, we present DEFT-FUNNEL, an open-source global optimization algorithm for general constrained grey-box and black-box problems that belongs to the class of trust-region sequential quadratic optimization algorithms. It extends the previous works by Sampaio and Toint (2015, 2016) to a global optimization solver that is able to exploit information from closed-form functions. Polynomial interpolation models are used as surrogates for the black-box functions and a clustering-based multistart strategy is applied for searching for the global minima. Numerical experiments show that DEFT-FUNNEL compares favorably with other state-of-the-art methods on two sets of benchmark problems: one set containing problems where every function is a black box and another set with problems where some of the functions and their derivatives are known to the solver. The code as well as the test sets used for experiments are available at the Github repository \href{http://github.com/phrsampaio/deft-funnel}{http://github.com/phrsampaio/deft-funnel}.
\end{abstract}


\section{Introduction}\label{sec:intro}

We are interested in finding the global minimum of the optimization problem
\begin{eqnarray}\label{eq:originalproblem}
\left\{\begin{array}{rc}
\displaystyle\min_x & f(x)\\
          \textrm{s.t.:} & l^c \leq c(x) \leq u^c, \\
          				 & l^h \leq h(x) \leq u^h, \\
          				 & l^x \leq x \leq u^x, \\
\end{array} \right.
\end{eqnarray}
where $f: \mathbb{R}^n \rightarrow \Re$ might be or not a black box, $c: \mathbb{R}^n \rightarrow
\mathbb{R}^q$ are black-box constraint functions and $h: \mathbb{R}^n \rightarrow \mathbb{R}^l$ are white-box constraint functions. i.e. their analytical expressions as well as their derivatives are available. The vectors $l^c$, $l^h$, $u^c$ and $u^h$ are lower and upper bounds on the constraints values $c(x)$ and $h(x)$, while $l^x$ and $u^x$ are bounds on the $x$ variables, with $l^c \in (\mathbb{R} \cup -\infty)^q$, $l^h \in (\mathbb{R} \cup -\infty)^l$, $u^c \in (\mathbb{R} \cup \infty)^q$, $u^h \in (\mathbb{R} \cup \infty)^l$, $l^x \in (\mathbb{R} \cup -\infty)^n$ and $u^x \in (\mathbb{R} \cup \infty)^n$. We also address the case where there are no white-box constraint functions $h$. Finally, we assume that the bound constraints are unrelaxable, i.e. feasibility must be maintained throughout the iterations, while the other general constraints are relaxable.

When at least one of the functions in (1) has a closed form, i.e. either the objective function is a white box or there is at least one white-box constraint function in the problem, the latter is said to be a grey-box problem. If no information about the functions is given at all, which means that the objective function is a black box and that there are no white-box constraints, the problem is known as a black-box problem. Both grey-box and black-box optimization belong to the field of derivative-free optimization (DFO) \cite{ConnScheinbergVincente09b,AudetHare2017}, where the derivatives of the functions are not available. DFO problems are encountered in real-life applications in various fields such as engineering, medicine, science and artificial intelligence. The black boxes are often the result of an expensive simulation or a proprietary code, in which case automatic differentiation \cite{Griewank2003,Griewank2008} is not applicable.

Many optimizations methods have been developed for finding stationary points or local minima of (1) when both the objective and the constraints functions are black boxes (e.g., \cite{Powell1994,Lewis2002,Bueno2013,Audet2015,Sampaio2015,Sampaio2016,Echebest2017,Amaioua2018}). However, a local minimum is not enough sometimes and so one need to search for the global minimum \cite{Floudas99,Floudas2000}. A reduced number of methods have been proposed to find the global minima of constrained black-box problems (see, for instance, \cite{Jones1998,Regis2005,Regis2007,Regis2011,Regis2014,Boukouvala2017}) and thus there is still much research to be done on this area.  Moreover, many global optimization methods for constrained black-box problems proposed in the literature or used in industry are unavailable to the public and are not open source. We refer the reader to the survey papers \cite{Rios2013,Larson2019} and to the textbooks \cite{ConnScheinbergVincente09b,AudetHare2017} for a comprehensive review on DFO algorithms for different types of problems.

In the case of constrained grey-box problems, especially those found in industry, it is a good idea to exploit the available information about the white boxes since such problems are usually hard to be solved, i.e. highly nonlinear, multimodal and with very expensive functions. Therefore, one would expect the solver to use any information given as input in order to attain the global minimum as fast as possible. Unfortunately, even less global optimization solvers exist for such problems today. A common solution found by optimization researchers, engineers and practioners is to consider all the functions as black boxes and to use a black-box optimization algorithm to solve the problem. Two of the few methods that exploit the available information are ARGONAUT \cite{Boukouvala2017} and a trust-region two-phase algorithm proposed in \cite{Bajaj2018}. In ARGONAUT, the black-box functions are replaced by surrogate models and a global optimization algorithm is used to solve the problem to global optimality, having the surrogate models updated only after its resolution. After updating the models, the problem is solved again with the updated models and this process is repeated until convergence is declared. In the two-phase algorithm described in \cite{Bajaj2018}, radial basis functions (RBF) are used as surrogate models for the black-box functions. Moreover, as in DEFT-FUNNEL, the self-correcting geometry approach proposed by \cite{ScheinbergToint10} is applied to the management of the interpolation set. Their algorithm is composed of a feasibility phase, where the goal is to find a feasible point, and an optimization phase, where the feasible point found in the first phase is used as a starting point to find a global minimum.

The algorithm in \cite{Bajaj2018} is the one sharing more elements in common with DEFT-FUNNEL. However, these two methods are still very different in nature since DEFT-FUNNEL combines a multistart strategy with a sequential quadratic optimization (SQO) algorithm in order to find global minima while the other applies a global optimization solver in its both phases for this purpose. Besides, DEFT-FUNNEL employs local polynomial interpolation models rather than RBF models as the former have good performance in the context of local optimization, which suits well for the SQO algorithm used in its local search. Despite the good performance of the algorithms proposed in \cite{Boukouvala2017} and \cite{Bajaj2018}, neither is freely available or open source.

\textbf{Contributions.} This paper proposes a new global optimization solver for general constrained grey-box and black-box problems written in Matlab \cite{Matlab2015} that exploits any available information given as input and that employs surrogate models built from polynomial interpolation in a trust-region-based SQO algorithm. Differently from the ARGONAUT approach, the surrogate models are updated during the optimization process as soon as new information from the evaluation of the functions at the iterates become available. Furthermore, the proposed solver, named DEFT-FUNNEL, is open source and freely available at the Github repository \href{http://github.com/phrsampaio/deft-funnel}{http://github.com/phrsampaio/deft-funnel}. It is based on the previous works in \cite{Sampaio2015, Sampaio2016} and it extends the original DFO algorithm to grey-box problems and to the search for global minima. As its previous versions, it does not require feasible starting points. To our knowledge, DEFT-FUNNEL is the first open-source global optimization solver for general constrained grey-box and black-box problems that exploits the derivative information available from the white-box functions. It is also the first one of the class of trust-funnel algorithms \cite{Gould2010} to be used in the search for global minima in both derivative-based and derivative-free optimization.

This paper serves also as the first release of the DEFT-FUNNEL code to the open-source community and to the public in general. For this reason and also due to the new extensions mentioned above, we give a complete description of the algorithm so that the reader can understand the method more easily while examining the code on Github. We also notice that some modifications and additions have been made to the local-search SQO algorithm with respect to the one presented in \cite{Sampaio2016}. In particular,  some changes were done in the condition for the normal step calculation, in the criticality step and in the maintenance of the interpolation set, all of them being described in due course. Furthermore, we have also added a second-order correction step.

The extension to global optimization is based on the multi-level single linkage (MLSL) method \cite{RinnooyKan1987a, RinnooyKan1987b}, a well-known stochastic multistart strategy that combines random sampling, a clustering-based approach for the selection of the starting points and local searches in order to identify all local minima. In DEFT-FUNNEL, it is used for selecting the starting points of the local searches done with the trust-funnel SQO algorithm.

\textbf{Organization.} The outline of this paper is as follows. Section \ref{sec:mlsl} introduces the MLSL method while in Section \ref{sec:thesolver} the DEFT-FUNNEL solver is presented in detail. In Section \ref{sec:numericalexp}, numerical results on a set of benchmark problems for global optimization are shown and the performance of DEFT-FUNNEL is compared with those of other state-of-the-art algorithms in a black-box setting. Moreover, numerical results on a set of grey-box problems are also analyzed. Finally, some conclusions about the proposed solver are drawn in Section \ref{sec:conclusions}.

\textbf{Notation.} Unless otherwise specified, the norm $\| \cdot \|$ is the standard Euclidean norm. Given any vector $x\in\mathbb{R}^n$, we denote its $i$-th component by $\left[x\right]_i$. We define $\left[x\right]^+=\max(0,x)$ where the $\max$ operation is done componentwise. We let $\calB(z; \Delta)$ denote the closed Euclidian ball centered at z, with radius $\Delta > 0$. Given any set $\calA$, $|\calA|$ denotes the cardinality of $\calA$. By $\calP_n^d$, we mean the space of all polynomials of degree at most $d$ in $\Re^n$. Finally, given any subspace $\calS$, we denote its dimension by $dim(\calS)$.

\section{Multi-level single linkage}\label{sec:mlsl}

The MLSL method \cite{RinnooyKan1987a, RinnooyKan1987b} is a stochastic multistart strategy originally designed for bound-constrained global optimization problems as below
\begin{eqnarray}
\left\{\begin{array}{rc}
\displaystyle\min & f(x)\\
          \textrm{s.t.:} & x\in\Omega, \\
\end{array} \right.
\end{eqnarray}
where $\Omega \subseteq \mathbb{R}^n$ is a convex, compact set containing the global minimum as an interior point that is defined by lower and upper bounds. It was later extended to problems with general constraints in \cite{Sendin2009}. As most of the stochastic multistart methods, it consists of a global phase, where random points are sampled from a probabilistic distribution, and a local phase, where selected points from the global phase are used as starting points for local searches done by some local optimization algorithm. MLSL aims at avoiding unnecessary and costly local searches that culminate in the same local minima. To achieve this goal, sample points are drawn from an uniform distribution in the global phase and then a local search procedure is applied to each of them except if there is another sample point or a previously detected local minimum within a critical distance with smaller objective function value. The method is fully described in Algorithm~\ref{MLSL}.

\algor{MLSL}{MLSL}{
%\vspace*{-0.3 cm}

\begin{description}
\item[1:] $\calL^{*} = \emptyset$
\item[2:] \textbf{for} $k=1$ to $\ldots$ \textbf{do}
\item[3:] \quad Generate $N$ random points uniformly distributed. 
\item[4:] \quad Rank the sample points by increasing value of $f$.
\item[5:] \quad \textbf{for} $i=1$ to $kN$ \textbf{do}
\item[6:] \quad\quad \textbf{if} $\not\exists x: \|x-x_i\| \leq r_k$ and $f(x) < f(x_i)$ \textbf{then}
\item[7:] \quad\quad\quad $\calL^{*} = \calL^{*} \cup LocalSearch(x_i)$
\item[8:] \quad\quad\textbf{end if}
\item[9:] \quad\textbf{end for}
\item[10:] \textbf{end for}
\item[11:] \textbf{return} the best local minimum found in $\calL^{*}$
\end{description} 
}
The critical distance $r_k$ is defined as
\begin{eqnarray}\label{eq:mlslradius}
r_k \eqdef \pi^{-1/2} \left(\Gamma \left(1 + \frac{n}{2}\right) m(S) \frac{\sigma \log kN}{kN} \right)^{1/n},
\end{eqnarray}
where $\Gamma$ is the gamma function and $m(S)$ is the Lebesgue measure of the set $S$. The method is centred on the idea of exploring the region of attraction of all local minima, which is formally defined below.

\theoremstyle{definition}
\begin{definition}\label{def:regionofattr}
Given a local search procedure $\calP$, we define a region of attraction $\calR(x^*)$ in $\Omega$ to be the set of all points in $\Omega$ starting from which $\calP$ will arrive at $x^*$.
\end{definition}

The ideal multistart method is the one that runs a local search only once at the region of attraction of every local minimum. However, two types of errors might occur in practice \cite{Locatelli1998}:
\begin{itemize}
    \item \textbf{Error 1.} The same local minimum $x^*$ has been found after applying local search to two or more points belonging to the same region of attraction of $x^*$.

  	\item \textbf{Error 2.} The region of attraction of a local minimum $x^*$ contains at least one sampled point, but local search has never been applied to points in this region.
\end{itemize}
In \cite{RinnooyKan1987a, RinnooyKan1987b}, the authors demonstrate the following theoretical properties of MLSL that are directly linked to the errors above:
\begin{itemize}
    \item \textbf{Property 1}. (Theorem 8 in \cite{RinnooyKan1987a} and Theorem 1 in \cite{RinnooyKan1987b}) If $\sigma>4$ in~\eqref{eq:mlslradius}, then, even if the sampling continues forever, the total number of local searches ever started by MLSL is finite with probability 1.

  	\item \textbf{Property 2}. (Theorem 12 in \cite{RinnooyKan1987a} and Theorem 2 in \cite{RinnooyKan1987b}) If $r_k$ tends to 0 with increasing $k$, then any local minimum $x^*$  will be found within a finite number of iterations with probability 1.
\end{itemize}
Property 1 states that the number of possible occurrences of Error 1 is finite while Property 2 says that Error 2 never happens. Due to its strong theoretical results and good practical performance, MLSL became one of the most reliable and popular multistart methods of late.

Finally, we note that MLSL has already been applied into the global constrained black-box optimization setting in \cite{Armstrong2016}. The proposed method, MLSL-MADS, integrates MLSL with a mesh adaptive direct search (MADS) method \cite{Audet2006} to find multiple local minima of an inverse transport problem involving black-box functions.

\section{The DEFT-FUNNEL solver}\label{sec:thesolver}

We first define the function $z: \mathbb{R}^n \rightarrow \mathbb{R}^{q+l}$ as $z(x) \eqdef (c(x), h(x))$, i.e. it includes all the constraint functions of the original problem \eqref{eq:originalproblem}. Then, by defining $f(x,s) \eqdef f(x)$ and $z(x,s) \eqdef z(x) - s$, the problem \eqref{eq:originalproblem} is rewritten as
\begin{eqnarray}
\left\{\begin{array}{rc}\label{eq:theproblem}
\displaystyle\min_{(x,s)} & f(x,s)\\
          \textrm{s.t.:} & z(x,s) = 0, \\
          				 & l^s \leq s \leq u^s, \\
          				 & l^x \leq x \leq u^x, \\
\end{array} \right.
\end{eqnarray}
where $s\in\mathbb{R}^{q+l}$ are slack variables and $l^s\in (\mathbb{R} \cup -\infty)^{q+l}$ and $u^s\in (\mathbb{R} \cup \infty)^{q+l}$ are the lower and upper bounds of the modified problem with $l^s = \left[l^c \; l^h\right]^T$ and $u^s = \left[u^c \; u^h\right]^T$. We highlight that the rewriting of the original problem \eqref{eq:originalproblem} as \eqref{eq:theproblem} is done within the solver and that the user does not need to interfere.

DEFT-FUNNEL is composed of a global search and a local search that are combined to solve the problem \eqref{eq:theproblem}. In the next two sections, we elaborate on each of these search steps.

\subsection{Global search}\label{subsec:globalsearch}

As mentioned previously, the global search in DEFT-FUNNEL relies on the MLSL multistart strategy. A merit function $\Phi(x)$ is used to decide which starting points are selected for the local search. We introduce the global search in detail in what follows.

\algor{GlobalSearch}{GlobalSearch}{
%\vspace*{-0.3 cm}

\begin{description}
\item[1:] $\calL^{*} = \emptyset$
\item[2:] \textbf{for} $k=1$ to $\ldots$ \textbf{do}
\item[3:] \quad Generate $N$ random points uniformly distributed. 
\item[4:] \quad Rank the sample points by increasing value of $\Phi$.
\item[5:] \quad \textbf{for} $i=1$ to $kN$ \textbf{do}
\item[6:] \quad\quad \textbf{if} $\not\exists x: \|x-x_i\| \leq r_k$ and $\Phi(x) < \Phi(x_i)$ \textbf{then}
\item[7:] \quad\quad\quad $\calL^{*} = \calL^{*} \cup LocalSearch(x_i)$
\item[8:] \quad\quad\textbf{end if}
\item[9:] \quad\textbf{end for}
\item[10:] \textbf{end for}
\item[11:] \textbf{return} the best feasible local minimum found in $\calL^{*}$
\end{description} 
}

The Algorithm \ref{GlobalSearch} is implemented in the function \verb|deft_funnel_multistart|, which is called by typing the following line in the Matlab command window:
\begin{lstlisting}
[best_sol, best_fval, best_indicators, total_eval, nb_local_searches, fL] = deft_funnel_multistart(@f, @c, @h, @dev_f, @dev_h, n, nb_cons_c, nb_cons_h)
\end{lstlisting}
The inputs and outputs of \verb|deft_funnel_multistart| are detailed in Table \ref{table:gsinout}. The options for the input `whichmodel' are described in Section \ref{subsec:localsearch}.

\begin{table}[ht]
\caption{Inputs and outputs of the function \texttt{deft\_funnel\_multistart}.}
\centering
\begin{tabular}{l l l}
\hline
                         & \textbf{Name} & \textbf{Description}\\
\hline
\textbf{Mandatory} & f & function handle of the objective function \\
\textbf{Inputs}    & c & function handle of the black-box constraints if any or an empty array \\
                     & h & function handle of the white-box constraints if any or an empty array \\
                     & dev\_f & function handle of the derivatives of f if it is a white box or an empty array \\
                     & dev\_h & function handle of the derivatives of h if any or an empty array \\
                     &  n & number of decision variables \\
                     & nb\_cons\_c & number of black-box constraints (bound constraints not included) \\
                     & nb\_cons\_h & number of white-box constraints (bound constraints not included) \\
\hline
\textbf{Optional}    & lsbounds & vector of lower bounds for the constraints \\
\textbf{Inputs}      & usbounds & vector of upper bounds for the constraints \\
                     & lxbounds & vector of lower bounds for the x variables \\
                     & uxbounds & vector of upper bounds for the x variables \\
                     & maxeval  & maximum number of evaluations (default: 5000*n) \\
                     & maxeval\_ls & maximum number of evaluations per local search (default: maxeval*0.7) \\
                     & whichmodel & approach to build the surrogate models \\
                     & f\_global\_optimum & known objective function value of the global optimum \\
\hline
\textbf{Outputs}     & best\_sol & best feasible solution found \\
                     & best\_fval & objective function value of ``best\_sol'' \\
                     & best\_indicators & indicators of ``best\_sol'' \\
                     & total\_eval & number of evaluations used \\
                     & nb\_local\_searches & number of local searches done \\
                     & fL & objective function values of all local minima found \\
\hline
\end{tabular}
\label{table:gsinout}
\end{table}

The merit function $\Phi(x)$ employed in DEFT-FUNNEL is the well-known $\ell_1$ penalty functon which is defined as follows:
\begin{eqnarray}\label{eq:merifunction}
\Phi(x) \eqdef f(x) + \pi\displaystyle\sum_{i=1}^{m}\left(\left[z(x) - u^s\right]^+ + \left[l^s - z(x)\right]^+\right),
\end{eqnarray}
where $m = q+l$ is the total number of constraints and $\pi$ is the penalty parameter. One of the advantages of this penalty function over others is that it is exact, that is, for sufficiently large values of $\pi$, the local minimum of $\Phi(x)$ subject only to the bound constraints on the $x$ variables is also the local minimum of the original constrained problem \eqref{eq:originalproblem}. We note that, although $\Phi(x)$ is nondifferentiable, it is only used in the global search for selecting the starting points for the local searches.

The \textit{LocalSearch} algorithm at line 7 is started by calling the function \verb|deft_funnel| whose inputs and outputs are given in the next section.

\subsection{Local search}\label{subsec:localsearch}

\begin{figure}[ht]
     \centering
     \includegraphics[scale=1]{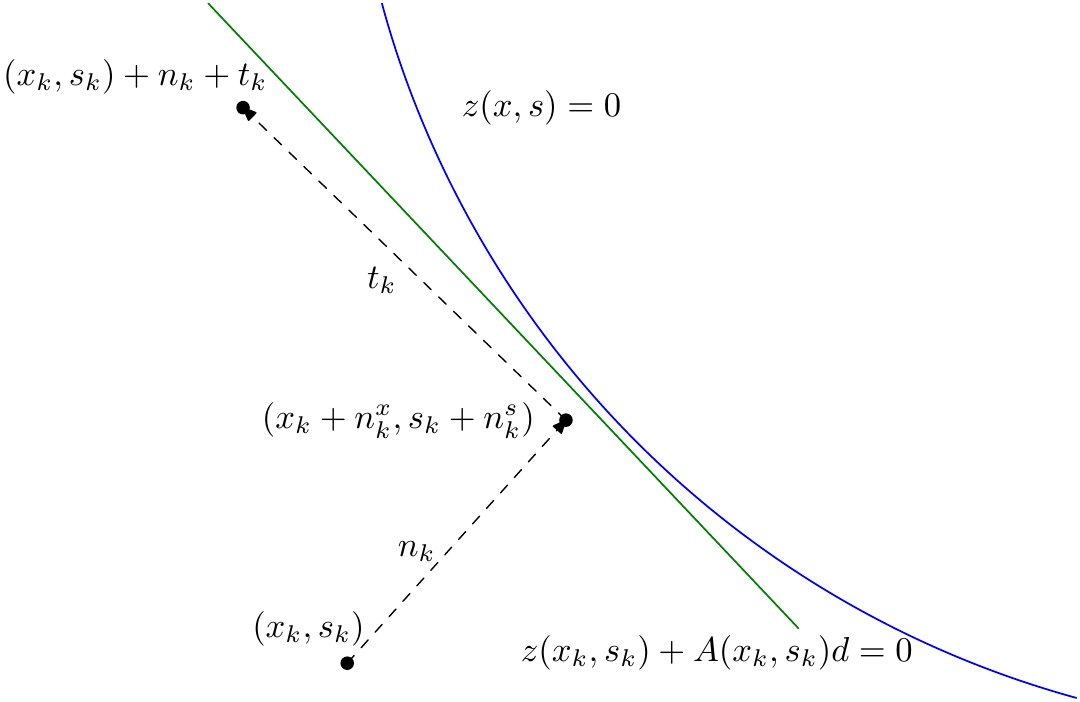}
     \caption{Illustration of the composite step $d_k = n_k + t_k$. The normal step $n_k$ attempts to improve feasibility by reducing the linearized constraint violation at $(x_k, s_k)$, whereas the tangent step aims at minimizing the objective function without deteriorating the gains in feasibility obtained through the normal step. Here, $A(x_k,s_k)$ denotes the Jacobian of $z(x,s)$ at $(x_k,s_k)$.}
     \label{fig:compositestep}
\end{figure}

The algorithm used in the local search is based on the one described in the paper \cite{Sampaio2016}. It is a trust-region SQO method that makes use of a funnel bound on the infeasibility of the iterates in order to ensure convergence. At each iteration $k$, we have an iterate $(x_k,s_k)$ such that $x_k\in\calY_k$, where $\calY_k$ defines the interpolation set at iteration $k$. Every iterate satisfies the following bound constraints:
\begin{eqnarray}
l^s \leq s_k \leq u^s, \label{eq:boundsslacks}\\
l^x \leq x_k \leq u^x. \label{eq:boundsx}
\end{eqnarray}
Depending on the optimality and feasibility of $(x_k,s_k)$, a new step $d_k \eqdef (d^x_k,d^s_k)^T$ is computed. Each full step of the trust-funnel algorithm is decomposed as 
\begin{eqnarray}\label{eq:stepdecomposition}
d_k = \left( \begin{array}{c}
d^x_k \\
d^s_k \end{array} \right)
= 
\left( \begin{array}{c}
n^x_k \\
n^s_k \end{array} \right)
+
\left( \begin{array}{c}
t^x_k \\
t^s_k \end{array} \right) = n_k + t_k,
\end{eqnarray}
where the normal step component $n_k$ aims to improve feasibility and the tangent step component $t_k$ reduces the objective function model without worsening the constraint violation up to first order. This is done by requiring the tangent step to lie in the null space of the Jacobian of the constraints and by requiring the predicted improvement in the objective function obtained in the tangent step to not be negligible compared to the predicted change in $f$ resulting from the normal step. The full composite step $d_k$ is illustrated in Figure~\ref{fig:compositestep}. As it is explained in the next subsections, the computation of the composite step in our algorithm does not involve the function $z(x.s)$ itself but rather its surrogate model.

After having computed a trial point $x_k + d_k$, the algorithm proceeds by checking whether the iteration was successful in a sense yet to be defined. The iterate is then updated correspondingly, while the sample set and the trust regions are updated according to a self-correcting geometry scheme to be described later on. If $\calY_k$ has been modified, the surrogate models are updated to satisfy the interpolation conditions for the new set $\calY_{k+1}$, implying that new function evaluations are carried out for the additional point obtained at iteration $k$.

The local search is started inside the multistart strategy loop in Algorithm \ref{GlobalSearch}, but it can also be called directly by the user in order to use DEFT-FUNNEL without multistart. This is done by typing the following line at the Matlab command window:
\begin{lstlisting}
[x, fx, mu, indicators, evaluations, iterate, exit_algo] = deft_funnel(@f, @c, @h, @dev_f, @dev_h, x0, nb_cons_c, nb_cons_h)
\end{lstlisting}
The inputs and outputs of the function \verb|deft_funnel| are detailed below in Table \ref{table:lsinout}. Many other additional parameters can be set directly in the function \verb|deft_funnel_set_parameters|. Those are related to the trust-region mechanism, to the interpolation set maintenance and to criticality step thresholds.

\begin{table}[ht]
\caption{Inputs and outputs of the function \texttt{deft\_funnel}}
\centering
\begin{tabular}{l l l}
\hline
                         & \textbf{Name} & \textbf{Description}\\
\hline
\textbf{Mandatory}   & f & function handle of the objective function \\
\textbf{Inputs}      & c & function handle of the black-box constraints if any or an empty array \\
                     & h & function handle of the white-box constraints if any or an empty array \\
                     & dev\_f & function handle of the derivatives of f if it is a white box or an empty array \\
                     & dev\_h & function handle of the derivatives of h if any or an empty array \\
                     & x0 & starting point (no need to be feasible) \\
                     & nb\_cons\_c & number of black-box constraints (bound constraints not included) \\
                     & nb\_cons\_h & number of white-box constraints (bound constraints not included) \\
\hline
\textbf{Optional}    & lsbounds & vector of lower bounds for the constraints \\
\textbf{Inputs}      & usbounds   & vector of upper bounds for the constraints \\
                     & lxbounds   & vector of lower bounds for the x variables \\
                     & uxbounds   & vector of upper bounds for the x variables \\
                     & maxeval    & maximum number of evaluations (default: 500*n) \\
                     & type\_f    & string 'BB' if f is a black box (default) or 'WB' otherwise \\
                     & whichmodel & approach to build the surrogate models \\
\hline
\textbf{Outputs}     & x & the best approximation found to a local minimum \\
                     & fx & the value of the objective function at x \\
                     & mu & local estimates for the Lagrange multipliers \\
                     & indicators & feasibility and optimality indicators \\
                     & evaluations & number of calls to the objective function and constraints \\
                     & iterate & info related to the best point found as well as the coordinates of all past iterates \\
                     & exit\_algo & output signal (0: terminated with success; -1: terminated with errors) \\
\hline
\end{tabular}
\label{table:lsinout}
\end{table}

Once the initial interpolation set has been built using one of the methods described in the next subsection, the algorithm calls the function \verb|deft_funnel_main|. In fact, \verb|deft_funnel| serves only as a wrapper for the main function of the local search, \verb|deft_funnel_main|, doing all the data preprocessing and paramaters setting needed in the initialization process. All the main steps such as the subspace minimization step, the criticality step and the computation of the new directions are part of the scope of \verb|deft_funnel_main|.

%%%%%%%%%%%%%%%%%%%%%%%%%%%%%%%%%%%%%%%%%%%%%%%%%%%%%%%%%%%%%%%%%%%%%%%%%%%%%%%%%%%%%%%%%

\subsubsection{Building the surrogate models}\label{subsubsec:surr}

The local-search algorithm starts by building an initial interpolation set either from a simplex or by drawing samples from an uniform distribution. The construction of the interpolation set is done in the function \verb|deft_funnel_build_initial_sample_set|, which is called only once during a local search within  \verb|deft_funnel|. The choice between random sampling and simplex is done in \verb|deft_funnel| when calling the function \verb|deft_funnel_set_parameters|, which defines the majority of parameters of the local search. If random sampling is chosen, it checks if the resulting interpolation set is well poised and, if not, it is updated using the Algorithm 6.3 described in Chapter 6 in \cite{ConnScheinbergVincente09b}, which is implemented in \verb|deft_funnel_repair_Y|.

For the sake of simplicity, we assume henceforth that the objective function is also a black box. Given a poised set of sample points $\calY_0=\{y^0,y^1,\ldots,y^p\}$ with an initial point $x_0\in \calY_0$, the next step of our algorithm is to replace the objective function $f(x)$ and the black-box constraint functions $c(x)=(c_1(x),c_2(x),\ldots,c_q(x))$ by surrogate models $m^f(x)$ and $m^c(x)=(m^{c_1}(x),m^{c_2}(x),\ldots,m^{c_q}(x))$ built from the solution of the interpolation system
\begin{equation}\label{eq:interpsystem}
M(\phi, \calY)\alpha_{\phi} = \Upsilon(\calY),
\end{equation}
where
\[
M(\phi,\calY) = \begin{pmatrix}
  \phi_{0}(y^0) & \phi_{1}(y^0) & \cdots & \phi_{b}(y^0) \\
  \phi_{0}(y^1) & \phi_{1}(y^1) & \cdots & \phi_{b}(y^1) \\
  \vdots  & \vdots  & \ddots & \vdots  \\
  \phi_{0}(y^p) & \phi_{1}(y^p) & \cdots & \phi_{b}(y^p)
 \end{pmatrix},
\;\;\Upsilon(\calY) = \begin{pmatrix}
  \Upsilon(y^0) \\
  \Upsilon(y^1) \\
  \vdots  \\
  \Upsilon(y^p)
 \end{pmatrix},
\;\;p \leq b,
\]
where $\phi$ is the basis of monomials and $\Upsilon(x)$ is replaced by the objective function $f(x)$ or some black-box constraint function $c_j(x)$.

We consider underdetermined quadratic interpolation models that are fully linear and that are enhanced with curvature information along the optimization process. Since the linear system \eqref{eq:interpsystem} is potentially underdetermined, the resulting interpolating polynomials will be no longer unique and so we provide to the user four approaches to construct the models $m^f(x)$ and $m^{c_j}(x)$ that can be chosen by passing a number from 1 to 4 to the input `whichmodel': 1 - subbasis selection approach; 2 - minimum $\ell_2$-norm models (by defatul); 3 - minimum Frobenius norm models; and 4 - regression mnodels (recommended for noisy functions). Further details about each one can be found in \cite{ConnScheinbergVincente09b}.

If $p_0=|\calY_0|=n+1$, a linear model  rather than an underdetermined quadratic model  is built for each function. The reason is that, despite both having error bounds that are linear in $\Delta$ for the first derivatives, the error bound for the latter includes also the norm of the model Hessian, as stated in Lemma 2.2 in \cite{Zhang2010}, which makes it worse than the former. 

Whenever $n+1 < p_k \leq (n+1)(n+2)/2 = p^{\max}$, the algorithm builds underdetermined quadratic models based on the choice of the user between the approaches described above. If regression models are considered instead, we set $p^{\max} = (n+1)(n+2)$, which means that the sample set is allowed to have twice the number of sample points required for fully quadratic interpolation models. Notice that having a number of sample points much larger than the required for quadratic interpolation can also worsen the quality of the interpolation models as the sample set could contain points that are too far from the iterate, which is not ideal for models built for local approximation.

It is also possible to choose the initial degree of the models 
between fully linear, quadratic with a diagonal Hessian or fully quadratic. This is done within \verb|deft_funnel| by setting the input argument \verb|cur_degree| in the call to \verb|deft_funnel_set_parameters| to one of the following options: \verb|model_size.plin|, \verb|model_size.pdiag| or \verb|model_size.pquad|.

The interpolation system \eqref{eq:interpsystem} is solved using a QR factorization of the matrix $M(\phi, \calY)$ within the function \verb|deft_funnel_computeP|, which is called by \verb|deft_funnel_build_models|.

In order to evaluate the error of the interpolation models and their derivatives with respect to the original functions $f$ and $c$, we make use of the measure of well poisedness of $\calY$ given below.

\theoremstyle{definition}
\begin{definition}\label{def:lambdapoisedness}
Let $\calY=\{y^0,y^1,\ldots,y^p\}$ be a poised interpolation set and $\calP_n^d$ be a space of polynomials of degree less than or equal to $d$ on $\mathbb{R}^n$. Let $\Lambda > 0$ and $\{\ell_0(x),\ell_1(x),\ldots,\ell_p(x)\}$ be the basis of Lagrange polynomials associated with $\calY$. Then, the set $\calY$ is said to be $\Lambda$-poised in $\calB$ for $\calP_n^d$ (in the interpolation sense) if and only if
\begin{eqnarray}
\displaystyle\max_{0\leq i \leq p}\displaystyle\max_{x\in\calB} |\ell_i(x)| \leq \Lambda. \nonumber
\end{eqnarray}
\end{definition}

As it is shown in \cite{ConnScheinbergVincente09b}, the error bounds for at most fully quadratic models depend linearly on the constant $\Lambda$; the smaller it is, the better the interpolation models approximate the original functions. We also note that the error bounds for undetermined quadratic interpolation models are linear in the diameter of the smallest ball containing $\calY$ for the first derivatives and quadratic for the function values. 

Finally, since the constraint functions $c(x)$ are replaced by surrogate models $m^c(x)$ in the algorithm, we define the models $m^z(x) \eqdef (m^c(x),h(x))$ and $m^z(x,s) \eqdef m^z(x) - s$, which are those used for computing new directions.

%%%%%%%%%%%%%%%%%%%%%%%%%%%%%%%%%%%%%%%%%%%%%%%%%%%%%%%%%%%%%%%%%%%%%%%%%%%%%%%%%%%%%%%%%

\subsubsection{Subspace minimization}

In this subsection, we explain how the subspace minimization is employed in our algorithm. We define the subspace $\calS_k$ at iteration $k$ as
\[
\calS_k \eqdef \{x \in \mathbb{R}^n\;|\; \left[x\right]_{i} = \left[l^x\right]_{i} \;\textrm{ for }\; i\in\calL_k \;\textrm{ and }\; \left[x\right]_{i} = \left[u^x\right]_{i} \;\textrm{ for }\; i\in\calU_k\},
\]
where $\calL_k \eqdef \{ i\;|\; \left[x_{k}\right]_{i} - \left[l^x\right]_{i} \leq \epsilon_b\}$ and $\calU_k \eqdef \{ i\;|\; \left[u^x\right]_{i} - \left[x_{k}\right]_{i} \leq \epsilon_b\}$ define the index sets of (nearly) active variables at their bounds, for some small constant $\epsilon_b >0$. If $\calL_k\neq\emptyset$ or $\calU_k\neq\emptyset$, a new well-poised interpolation set $\calZ_k$ is built and a recursive call is made in order to solve the problem in the new subspace $\calS_k$. 

\begin{figure}[ht]
     \centering
     \includegraphics[scale=0.9]{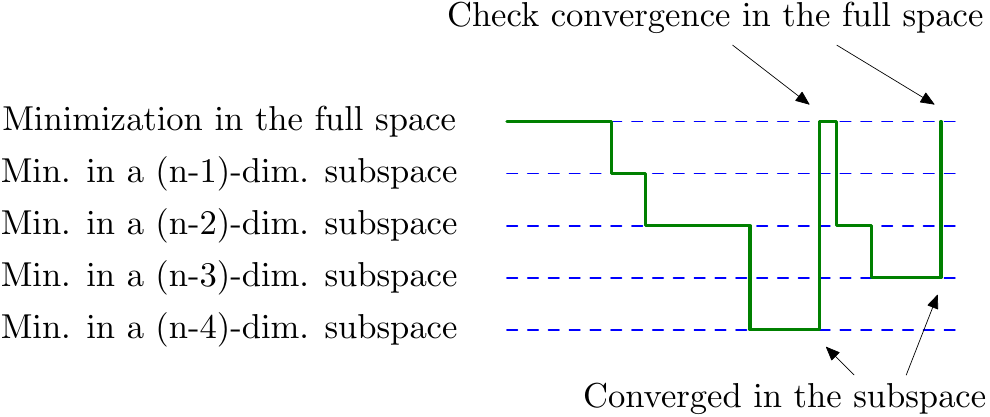}
     \caption{Subspace minimization procedure. The algorithm calls itself recursively in the order to solve the problem in a new subspace. If convergence is attained, it goes back to check if the solution found is also optimal in the full space.}
     \label{fig:subspacemin}
\end{figure}

If the algorithm converges in a subspace $\calS_k$ with an optimal solution $(\tilde{x},\tilde{s})$, it checks if the latter is also optimal for the full-space problem, in which case the algorithm stops. If not, the algorithm continues by attempting to compute a new direction in the full space. This procedure is illustrated in Figure~\ref{fig:subspacemin}.

\begin{figure}[ht]
     \centering
     \includegraphics[scale=1]{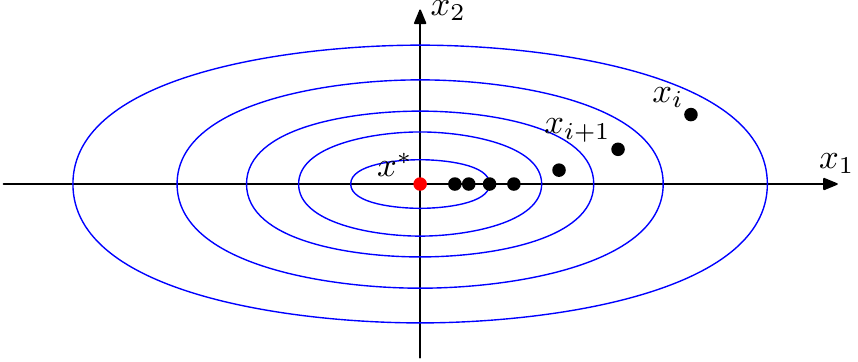}
     \caption{Illustration of a scenario where the interpolation set becomes degenerated as the optimal solution is approached. In this example, we consider a two-dimensional problem with the bound constraint $\left[x\right]_2\geq0$, which is active at the solution $x^*$ and at the iterates close to it.}
     \label{fig:degenset}
\end{figure}

The dimensionality reduction of the problem mitigates the chances of degeneration of the interpolation set when the sample points become too close to each other and thus affinely dependent. Figure~\ref{fig:degenset} gives an example of this scenario as the optimal solution is approached.

In order to check the criticality in the full-space problem, a full-space interpolation set of degree $n+1$ is built in an $\epsilon$-neighborhood around the point $x^*_{\calS}$, which is obtained by assembling the subspace solution $\tilde{x}$ and the $|\calL_k\,\cup\,\calU_k|$ fixed components $\left[x_k\right]_i$. The models $m^f_k$ and $m^c_k$ are then updated and the criticality step is entered.

The complete subspace minimization step is described in Algorithm~\ref{DFTRSUBSPACE} and it is implemented in the function \verb|deft_funnel_subspace_min|, which is called inside \verb|deft_funnel_main|.

\algor{DFTRSUBSPACE}{SubspaceMinimization($\calS_{k-1}$, $\calY_k$, $x_k$, $s_k$, 
$\Delta_k^f$, $\Delta_k^z$, $v_k^{\max}$)}{
%\vspace*{-0.3 cm}

\begin{description}
\item[1:] Check for (nearly) active bounds at $x_k$ and define $\calS_k$. If there is no (nearly) active bound or if $\calS_k$ has already been explored, go to Step 5. If all bounds are active, go to Step 4.
\item[2:] Build a new interpolation set $\calZ_k$ in $\calS_k$.
\item[3:] Call recursively \underline{LocalSearch}($\calS_k$, $\calZ_k$, $x_k$, $s_k$, 
$\Delta_k^f$, $\Delta_k^z$, $v_k^{\max}$) and let $(x^*_{\calS},s^*_{\calS})$ be the solution of the subspace problem after adding the fixed components.
\item[4:] If $dim(\calS_{k-1})<n$, return $(x^*_{\calS},s^*_{\calS})$. Otherwise, reset $(x_k,s_k)=(x^*_{\calS},s^*_{\calS})$, construct new set $\calY_k$ around $x_k$, build $m_k^f$ and $m_k^c$ and recompute $\pi^f_{k-1}$ (optimality measure).
\item[5:] If $\calS_k$ has already been explored, set $(x_{k+1},s_{k+1})=(x_k,s_k)$, reduce the trust regions radii $\Delta_{k+1}^f = \gamma\Delta_k^f$ and $\Delta_{k+1}^z=\gamma\Delta_k^z$, set $\Delta_{k+1}=\min[\Delta_{k+1}^f,\Delta_{k+1}^z]$ and build a new poised set $\calY_{k+1}$ in $\calB(x_{k+1};\Delta_{k+1})$.
\end{description} 
}

%%%%%%%%%%%%%%%%%%%%%%%%%%%%%%%%%%%%%%%%%%%%%%%%%%%%%%%%%%%%%%%%%%%%%%%%%%%%%%%%%%%%%%%%%

\subsubsection{The normal step}\label{subsubsec:normalstep}

The normal step aims at reducing the constraint violation at given point $(x,s)$ defined by
\begin{eqnarray}\label{eq:vdef}
v(x,s) \eqdef \half \| m^z(x,s) \|^2.
\end{eqnarray}
To ensure that the step $n_k$ is normal to the linearized constraint $m^z(x_k,s_k) + J(x_k,s_k) n = 0$, where $J(x,s)\eqdef (J(x)\;-I)$ is the Jacobian of $m^z(x,s)$ with respect to $(x,s)$ and $J(x)$ is the Jacobian of $m^z(x)$ with respect to $x$, we require that
\begin{eqnarray}\label{eq:nkbound}
\| n_k \|_{\infty} \leq \kappa_{n} \|m^z(x_k,s_k)\|,
\end{eqnarray}
for some $\kappa_{n} \geq 1$.

The computation of $n_k$ is done by solving the constrained linear least-squares subproblem
\begin{eqnarray}\label{eq:nstep}
\left\{\begin{array}{rc}
\displaystyle\min_{n=(n^x,n^s)} & \half\| m^z(x_k,s_k) + J(x_k,s_k) n\|^2 \\
          \textrm{s.t.:}  & l^s \leq s_k + n^s \leq u^s, \label{eq:nstepcon1}\\
                         & l^x \leq x_k + n^x \leq u^x,\label{eq:nstepcon2}\\
                         & x_k + n^x \in \calS_k, \label{eq:nstepcon3}\\
                         & n \in \calN_k, \label{eq:nstepcon4}
\end{array}\right .
\end{eqnarray}
where
\begin{eqnarray}\label{eq:ninN}
\calN_k \eqdef \{ n \in \mathbb{R}^{n+m} \mid \| n \|_{\infty} \leq \min\left[\,\Delta^z_k,\,\kappa_{n}\, \|m^z(x_k,s_k)\|\,\right]\, \},
\end{eqnarray}
for some trust-region radius $\Delta_k^z > 0$. Finally, a funnel bound $v_k^{\max}$ is imposed on the constraint violation $v_k \eqdef v(x_k,s_k)$ for the acceptance of new iterates to ensure the convergence towards feasibility.

We notice that, although a linear approximation of the constraints is used for calculating the normal step, the second-order information of the quadratic interpolation model $m^z(x,s)$ is still used in the SQO model employed in the tangent step problem as it is shown next.

The subproblem~\eqref{eq:nstep} is solved within the function \verb|deft_funnel_normal_step|, which makes use of an original active-set algorithm where the unconstrained problem is solved at each iteration in the subspace defined by the currently active bounds, themselves being determined by a projected Cauchy step. Each subspace solution is then computed using a SVD decomposition of the reduced matrix. This algorithm is implemented in the function \verb|deft_funnel_blls| and is intended for small-scale bound-constrained linear least-squares problems.

%%%%%%%%%%%%%%%%%%%%%%%%%%%%%%%%%%%%%%%%%%%%%%%%%%%%%%%%%%%%%%%%%%%%%%%%%%%%%%%%%%%%%%%%%

\subsubsection{The tangent step}\label{subsubsec:tangentstep}

The tangent step is a direction that improves optimality and it is computed by using a SQO model for the problem~\eqref{eq:theproblem} after the normal step calculation. The quadratic model for the function objective function is defined as
\begin{eqnarray}\label{eq:tstepmodel1}
\psi_k((x_k,s_k) + d) \eqdef m^f(x_k,s_k) + \ip{g_k}{d} + \half \ip{d}{B_k d},
\end{eqnarray}
where
$m^f(x_k,s_k) \eqdef m^f(x_k)$, $g_k \eqdef \nabla_{(x,s)} m^f(x_k,s_k)$, and $B_k$ is the approximate Hessian of the Lagrangian function
\begin{align}
\calL (x,s,\mu,\xi^s,\tau^s,\xi^x,\tau^x) & = m^f(x,s) + \ip{\mu}{m^z(x,s))} + \ip{\tau^s}{s - u^s} + \ip{\xi^s}{l^s - s} \nonumber \\
& \quad\; +  \ip{\tau^x}{x - u^x} + \ip{\xi^x}{l^x - x} \nonumber
\end{align}
with respect to $(x,s)$, given by 
\begin{eqnarray}\label{eq:defBk}
B_k = \begin{pmatrix}
  H_k + \sum_{i=1}^m [\hat{\mu}_k]_i Z_{ik} & 0 \\
  0   & 0
 \end{pmatrix},
\end{eqnarray}
where $\xi^s$ and $\tau^s$ are the Lagrange multipliers associated to the lower and upper bounds, respectively, on the slack variables $s$, and $\xi^x$ and $\tau^x$, the Lagrange multipliers associated to the lower and upper bounds on the $x$ variables. In \eqref{eq:defBk}, $H_k = \nabla^2_{xx} m^f(x_k,s_k)$, $Z_{ik} = \nabla^2_{xx} m^z_{ik}(x_k,s_k)$
and the vector $\hat{\mu}_k$ may be viewed as a local approximation of
the Lagrange multipliers with respect to the equality constraints $m^z(x,s) = 0$.

By applying \eqref{eq:stepdecomposition} into \eqref{eq:tstepmodel1}, we obtain
\begin{equation}\label{eq:tstepmodel3}
\begin{array}{rcl}
\psi_k((x_k,s_k) + n_k + t) 
& = & \psi_k((x_k,s_k) + n_k) + \ip{g_k^{N}}{t} + \half \ip{t}{B_k t},
\end{array}
\end{equation}
where
\begin{eqnarray}\label{eq:gkndef}
g_k^{N} \eqdef g_k + B_k\, n_k.
\end{eqnarray}
Since \eqref{eq:tstepmodel3} is a local approximation for the function $m^f((x_k,s_k)+n_k+t)$, a trust region with radius $\Delta_k^f$ is used for the complete step $d=n_k+t$:
\begin{eqnarray}\label{eq:trfdef}
\calT_k \eqdef \{ d \in \mathbb{R}^{n+m} \mid \| d \|_{\infty} \leq \Delta^f_k \}.
\end{eqnarray}
Moreover, given that the normal step was also calculated using local models, it makes sense to remain in the intersection of both trust regions, which implies that
\begin{eqnarray}\label{eq:sinB}
d_k \in \calR_k
\eqdef \calN_k \cap \calT_k
\eqdef \{ d \in \mathbb{R}^{n+m} \mid \| d \|_{\infty} \leq \Delta_k \},
\end{eqnarray}
where $\Delta_k = \min[ \Delta_k^z, \Delta_k^f]$.

In order to make sure that there is still enough space left for the tangent step within $\calR_k$, we first check if the following constraint on the normal step is satisfied: 
\begin{eqnarray}\label{eq:suffinside}
\|n_k\|_{\infty} \leq \kappa_{\calR} \Delta_k,
\end{eqnarray}
for some $\kappa_{\calR} \in (0,1)$. If~\eqref{eq:suffinside} holds, the tangent step is calculated by solving the following subproblem

\begin{eqnarray}\label{eq:tstep}
\left\{\begin{array}{rc}
\displaystyle\min_{t=(t^x,t^s)} & \ip{g_k^{N}}{t} + \half \ip{t}{B_k t} \\
\textrm{s.t.:}    & J(x_k,s_k) t = 0, \\
          	      & l^s \leq s_k + n^s_k + t^s \leq u^s, \\
                  & l^x \leq x_k + n^x_k + t^x \leq u^x, \\
                  & x_k + n^x_k + t^x \in \calS_k. \\
                  & n_k+t \in \calR_k, \\
\end{array} \right.
\end{eqnarray}
where we require that the new iterate $x_k+d_k$ belongs to subspace $\calS_k$ and that it satisfies the bound constraints ~\eqref{eq:boundsslacks} and~\eqref{eq:boundsx}. In the Matlab code, the tangent step is calculated by the function \verb|deft_funnel_tangent_step|, which in turn calls either the Matlab solver \verb|linprog| or our implementation of the nonmonotone spectral projected gradient method \cite{Birgin00} in \verb|deft_funnel_spg| to solve the subproblem \eqref{eq:tstep}. The choice between both solvers is based on whether $\|B_k\|\leq \epsilon$, for a small $\epsilon>0$, in which case we assume that the problem is linear and therefore \verb|linprog| is used.

We define our $f$-criticality measure as 
\begin{eqnarray}\label{eq:optmeasuredef}
\pi_k^f\eqdef -\ip{g_k^N}{r_k},
\end{eqnarray}
where $r_k$ is the projected Cauchy direction obtained by solving the linear optimization problem
\begin{eqnarray}\label{eq:optmeasureprob}
\left\{\begin{array}{rc}
\displaystyle\min_{r=(r^x,r^s)} & \ip{g_k^N}{r} \\
       \textrm{s.t.:} & J(x_k,s_k) r = 0, \\
          	          & l^s \leq s_k + n^s_k + r^s \leq u^s, \\
                      & l^x \leq x_k + n^x_k + r^x \leq u^x, \\
                      & x_k + n^x_k + r^x \in \calS_k. \\
                      & \|r\|_{\infty} \leq 1.
\end{array}\right .
\end{eqnarray}
By definition, $\pi_k^f$ measures how much decrease could be obtained locally
along the projection of the negative of the approximate gradient $g_k^N$ onto the nullspace of $J(x_k,s_k)$ intersected with the region delimited by the bounds. This measure is computed in \verb|deft_funnel_compute_optimality|, which uses \verb|lingprog| to solve the subproblem \eqref{eq:optmeasureprob}.

A new local estimate of the Lagrange
multipliers $(\mu_k,\xi^s_k,\tau^s_k,\xi^x_k,\tau^x_k)$ are computed by solving the following problem:
\begin{eqnarray}\label{eq:LSprob}
\left\{\begin{array}{rc}
\displaystyle\min_{(\mu,\hat{\xi}^s,\hat{\tau}^s,\hat{\xi}^x,\hat{\tau}^x)} & \half \| \calM_k(\mu,\hat{\xi}^s,\hat{\tau}^s,\hat{\xi}^x,\hat{\tau}^x) \|^2\\
\textrm{s.t.:}    & \hat{\xi}^s,\hat{\tau}^s,\hat{\xi}^x,\hat{\tau}^x \geq 0,
\end{array} \right.
\end{eqnarray}
where
\begin{eqnarray}\label{eq:lagmult}
\calM_k(\mu,\hat{\xi}^s,\hat{\tau}^s,\hat{\xi}^x,\hat{\tau}^x) & 
\eqdef & \left( \begin{array}{c}
g_k^{N} \\
0 \end{array} \right)
+
\left( \begin{array}{c}
J(x_k)^T \\
-I \end{array} \right)\mu
+
\left( \begin{array}{c}
0 \\
I_{\tau}^s \end{array} \right)\hat{\tau}^s
+
\left( \begin{array}{c}
0 \\
-I_{\xi}^s \end{array} \right)\hat{\xi}^s \nonumber \\
&      & +
\left( \begin{array}{c}
I_{\tau}^x \\
0 \end{array} \right)\hat{\tau}^x
+
\left( \begin{array}{c}
-I_{\xi}^x \\
0 \end{array} \right)\hat{\xi}^x, \nonumber
\end{eqnarray}
the matrix $I$ is the $m\times m$ identity matrix, the matrices $I_{\xi}^s$ and $I_{\tau}^s$ are obtained from $I$ by removing the columns whose indices are not associated to any active (lower and upper, respectively) bound at $s_k + n^s_k$, the matrices $I_{\xi}^x$ and $I_{\tau}^x$ are obtained from the $n\times n$ identity matrix by removing the columns whose indices are not associated to any active (lower and upper, respectively) bound at $x_k + n^x_k$, and the Lagrange multipliers $(\hat{\xi}^s,\hat{\tau}^s,\hat{\xi}^x,\hat{\tau}^x)$ are those in $(\xi^s,\tau^s,\xi^x,\tau^x)$ associated to active bounds at $s_k + n^s_k$ and $x_k + n^x_k$. All the other Lagrange multipliers are set to zero.

The subproblem~\eqref{eq:LSprob} is also solved using the active-set algorithm implemented in the function \verb|deft_funnel_blls|.

%%%%%%%%%%%%%%%%%%%%%%%%%%%%%%%%%%%%%%%%%%%%%%%%%%%%%%%%%%%%%%%%%%%%%%%%%%%%%%%%%%%%%%%%%

\subsubsection{Which steps to compute and retain}\label{subsubsec:whichsteps}
 
The algorithm computes normal and tangent steps depending on the measures of feasibility and optimality at each iteration. Differently from \cite{Sampaio2015,Sampaio2016}, where the computation of the normal steps depends on the measure of optimality, here the normal step is computed whenever the following condition holds
\begin{eqnarray}\label{eq:normal_cond}
\|z(x_k,s_k)\| > \epsilon,
\end{eqnarray}
for some small $\epsilon>0$, i.e. preference is always given to feasibility. This choice is based on the fact that, in many real-life problems with expensive functions and a small budget, one seeks to find a feasible solution as fast as possible and that a solution having a smaller objective function value than the current strategy is already enough. If \eqref{eq:normal_cond} fails, we set $n_k=0$. 

We define a $v$-criticality measure that
indicates how much decrease could be obtained locally
along the projection of the negative gradient of the Gauss-Newton model of $v$ at $(x_k,s_k)$ onto the region delimited by the bounds as 
\[
\pi_k^v \eqdef -\ip{J(x_k,s_k)^T z(x_k,s_k)}{b_k},
\]
 where the projected Cauchy direction $b_k$ is given by the solution of
\begin{eqnarray}\label{eq:vmeasureprob}
\left\{\begin{array}{rc}
\displaystyle\min_{b=(b^x,b^s)} & \ip{J(x_k,s_k)^T z(x_k,s_k)}{b} \\
       \textrm{s.t.:} & l^s \leq s_k + b^s \leq u^s, \\
                      & l^x \leq x_k + b^x \leq u^x, \\
                      & x_k + b^x \in \calS_k, \\
                      & \|b\|_{\infty} \leq 1.
\end{array}\right .
\end{eqnarray}
We say that $(x_k,s_k)$ is an infeasible stationary point if $z(x_k,s_k) \neq 0$ and $\pi_k^v=0$, in which case the algorithm terminates. 

The procedure for the calculation of the normal step is given in the algorithm below. In the code, it is implemented in the function \verb|deft_funnel_normal_step|, which calls the algorithm in \verb|deft_funnel_blls| in order to solve the normal step subproblem \eqref{eq:nstep}.

\algor{DFTRNORMALSTEP}{NormalStep($x_k$, $s_k$, $\pi_k^v$, $v_k$, $v_k^{\max}$)}{
%\vspace*{-0.3 cm}
\begin{description}
\item[1:] If $z(x_k,s_k) \neq 0$ and $\pi_k^v=0$, \textbf{STOP} (infeasible stationary point).
\item[2:] If \eqref{eq:normal_cond}, compute a normal step $n_k$ by solving \eqref{eq:nstep}. Otherwise, set $n_k=0$.
\end{description}
}

If the solution of \eqref{eq:optmeasureprob} is $r_k=0$, then by \eqref{eq:optmeasuredef} we have $\pi_k^f=0$, in which case we set $t_k=0$. If the current iterate is farther from feasibility than from optimality, i.e., for a given a monotonic bounding function $\omega_t$, the condition
\begin{eqnarray}\label{eq:fcond}
\pi_k^f > \omega_t(\|z(x_k,s_k)\|)
\end{eqnarray}
fails, then we skip the tangent step computation by setting $t_k=0$. 

After the computation of the tangent step, the usefulness of the latter is evaluated by checking if the following conditions
\begin{eqnarray}\label{eq:not-tkOnk}
\|t_k\| > \kappa_{\calZ\calS}\|n_k\|
\end{eqnarray}
and
\begin{eqnarray}\label{eq:fmodratios}
\delta_k^f
\eqdef \delta^{f,t}_k + \delta^{f,n}_k
\geq \kappa_\delta \delta^{f,t}_k,
\end{eqnarray}
where
\begin{eqnarray}\label{eq:tangent_Cauchy}
\delta^{f,t}_k
\eqdef \psi_k((x_k,s_k)+n_k) -\psi_k((x_k,s_k)+n_k+t_k)
\end{eqnarray}
and
\begin{eqnarray}\label{eq:delta_f_n}
\delta^{f,n}_k
\eqdef \psi_k(x_k,s_k) -\psi_k((x_k,s_k)+n_k),
\end{eqnarray}
are satisfied for some $\kappa_{\calZ\calS} > 1$ and for
$\kappa_\delta \in (0,1)$.  The inequality~\eqref{eq:fmodratios} indicates that the \emph{predicted}
improvement in the objective function obtained in the tangent step is substantial compared to the \emph{predicted} change in $f$ resulting from the normal step. If \eqref{eq:not-tkOnk} holds but \eqref{eq:fmodratios} fails, the tangent step is not useful in the sense just discussed, and we choose to ignore it by resetting $t_k=0$.

The tangent step procedure is stated in Algorithm~\ref{DFTRTANGSTEP} and it is implemented in the function \verb|deft_funnel_tangent_step|.

\algor{DFTRTANGSTEP}{TangentStep($x_k$, $s_k$, $n_k$)}{
%\vspace*{-0.3 cm}
\begin{description}
\item[1:] If \eqref{eq:suffinside} holds, then
\vspace*{-2mm}
\begin{description}
\item[1.1:] select a vector $\hat{\mu}_k$ and define $B_k$ as in \eqref{eq:defBk};
\item[1.2:] compute $\mu_k$ by solving \eqref{eq:LSprob}; 
\item[1.3:] compute the modified Cauchy direction $r_k$ by
  solving \eqref{eq:optmeasureprob} and define $\pi_k^f$ as \eqref{eq:optmeasuredef}; 
\item[1.4:]  if \eqref{eq:fcond} holds, compute a tangent step $t_k$ by solving \eqref{eq:tstep}.
\end{description}
\item[2:] If \eqref{eq:suffinside} fails, set $\mu_k=\mu_{k-1}$.
In this case, or if \eqref{eq:fcond} fails, or if \eqref{eq:not-tkOnk} holds but
\eqref{eq:fmodratios} fails, set $t_k=0$ and $d_k=n_k$.
\item[3:] Define $(x_k^+,s_k^+) = (x_k,s_k) + d_k$.
\end{description}
}

%%%%%%%%%%%%%%%%%%%%%%%%%%%%%%%%%%%%%%%%%%%%%%%%%%%%%%%%%%%%%%%%%%%%%%%%%%%%%%%%%%%%%%%%%

\subsubsection{Iteration types}\label{subsubsec:iterationstype}

Depending on the contributions of the current iteration in terms of optimality and feasibility, we classify it into three types: $\mu$-iteration, $f$-iteration and $z$-iteration. This is done by checking if some conditions hold for the trial point defined as
\begin{eqnarray}\label{eq:xk+def}
(x_k^+,s_k^+) \eqdef (x_k,s_k) + d_k.
\end{eqnarray}

\paragraph{$\mu$-iteration}

If $d_k=0$, which means that the Lagrange multiplier estimates are the only new values that have been computed, iteration $k$ is said to be a $\mu$-iteration with reference to the Lagrange multipliers $\mu$ associated to the constraints $m^z(x,s)=0$. Notice, however, that not only new $\mu_k$ values have been computed, but all the other Lagrange multipliers $(\xi^s_k,\tau^s_k,\xi^x_k,\tau^x_k)$ as well.

In this case, we set $(x_{k+1},s_{k+1}) = (x_k,s_k)$, $\Delta_{k+1}^f = \Delta_k^f$, $\Delta_{k+1}^z = \Delta_k^z$, $v_{k+1}^{\max} =
v_k^{\max}$ and we use the new multipliers to build a new SQO model in \eqref{eq:tstepmodel1}.

Since null steps $d_k=0$ might be due to the poor quality of the interpolation models, we check the $\Lambda$-poisedness in $\mu$-iterations and attempt to improve it whenever the following condition holds
\begin{equation}\label{eq:muiterpoisedness}
\Lambda\, \Delta(\calY_k)>\epsilon_{\mu},
\end{equation}
where 
\[
\Delta(\calY_k) \eqdef \displaystyle\max_j \|y^{k,j}-x_k\|
\]
and $\epsilon_{\mu}>0$. The inequality \eqref{eq:muiterpoisedness} gives an estimate of the error bound for the interpolation models. If \eqref{eq:muiterpoisedness} holds, we try to reduce the value at the left side by modifying the sample set $\calY_k$. Firstly, we choose a constant $\xi\in(0,1)$ and replace all points $y^{k,j}\in\calY_k$ such that 
\[
\|y^{k,j}-x_k\|>\xi\Delta(\calY_k)
\]
by new points $y^{k,j}_{*}$ that (approximately) maximizes $|\ell_{j_k}(x)|$ in $\calB(x_k;\xi\Delta(\calY_k))$. Then we use the Algorithm 6.3 described in Chapter 6 in \cite{ConnScheinbergVincente09b} with the smaller region $\calB$ to improve $\Lambda$-poisedness of the new sample set. This procedure is implemented in the Matlab code by the function \verb|deft_funnel_repair_sample_set|.

\paragraph{$f$-iteration}

If iteration $k$ has mainly contributed to optimality, we say that iteration $k$ is an $f$-iteration. Formally, this happens when $t_k \neq 0$, \eqref{eq:fmodratios} holds, and
\begin{eqnarray}\label{eq:feasok}
v(x_k^+,s_k^+) \leq v_k^{\max}.
\end{eqnarray}
Convergence of the algorithm towards feasibility is ensured by condition \eqref{eq:feasok}, which limits the constraint violation with the funnel bound.

In this case, we set $(x_{k+1},s_{k+1}) = (x_k^+,s_k^+)$) if
\begin{eqnarray}\label{eq:rhof}
\rho_k^f
\eqdef \frac{m^f(x_k,s_k) - m^f(x_k^+,s_k^+)}{\delta^f_k}
\geq \eta_1,
\end{eqnarray}
and $(x_{k+1},s_{k+1}) = (x_k,s_k)$) otherwise. Note that $\delta_k^f>0$ (because of \eqref{eq:tangent_Cauchy} and \eqref{eq:fmodratios}) unless $(x_k,s_k)$ is first-order critical, and hence that condition \eqref{eq:rhof} is well-defined. As for the value of the funnel bound, we set $v_{k+1}^{\max} = v_k^{\max}$. 

Since our method can suffer from the Maratos effect \cite{Maratos1978}, we also apply a second-order correction (see Chapter 15, Section 15.6, in \cite{Nocedal2006} for more details) to the normal step whenever the complete direction $d_k$ is unsuccessful at improving optimality, i.e., whenever the condition \eqref{eq:rhof} fails. As the latter might be due to the local approximation of the constraint functions, this effect may be overcome with a second-order step $\hat{n}$ that is calculated in \verb|deft_funnel_sec_order_correction| by solving the following subproblem
\begin{eqnarray}\label{eq:nstepsecorder}
\left\{\begin{array}{rc}
\displaystyle\min_{\hat{n}=(\hat{n}^x,\hat{n}^s)} & \half\| m^z(x_k^+,s_k^+) + J(x_k,s_k) \hat{n}\|^2 \\
          \textrm{s.t.:}  & l^s \leq s_k^+ + \hat{n}^s \leq u^s, \\
                         & l^x \leq x_k^+ + \hat{n}^x \leq u^x, \\
                         & x_k^+ + \hat{n}^x \in \calS_k, \\
                         & \hat{n} \in \hat{\calN_k},
\end{array}\right .
\end{eqnarray}
where
\begin{eqnarray}\label{eq:nhatinN}
\hat{\calN_k} \eqdef \{ \hat{n} \in \mathbb{R}^{n+m} \mid \| \hat{n} \|_{\infty} \leq \min\left[\,\Delta^z_k,\,\kappa_{n}\, \|m^z(x_k^+,s_k^+)\|\,\right]\, \}.
\end{eqnarray}

\paragraph{$z$-iteration}

If iteration $k$ is neither a $\mu$-iteration nor a $f$-iteration, then it is said to be a
$z$-iteration. This means that the major contribution of iteration $k$ is to improve feasibility, which happens when either $t_k=0$ or when \eqref{eq:fmodratios} fails.

We accept the trial point if the improvement in feasibility is comparable to its predicted value
\[
\delta^z_k \eqdef \half \|z(x_k,s_k)\|^2 - \half \| z(x_k,s_k) + J(x_k,s_k) d_k \|^2,
\]
and the latter is itself comparable to its predicted decrease along the normal step, that is
\begin{eqnarray}\label{eq:rhoc}
n_k \neq 0,\; \text{ } \delta_k^z \geq \kappa_{zn} \delta^{z,n}_k \tim{and} \rho_k^z \eqdef \frac{v(x_k,s_k) - v(x_k^+,s_k^+)}{\delta^z_k} \geq \eta_1,
\end{eqnarray}
for some $\kappa_{zn}\in (0,1-\kappa_{tg}]$ and where
\begin{eqnarray}\label{eq:normal_Cauchy}
\delta^{z,n}_k & 
\eqdef         & \half \|z(x_k,s_k)\|^2 - \half \|z(x_k,s_k)+ J(x_k,s_k) n_k\|^2.
\end{eqnarray}
If \eqref{eq:rhoc} fails, the trial point is rejected. 

Finally, the funnel bound is updated as follows
\begin{eqnarray}\label{eq:ccmax}
v_{k+1}^{\max} =
 \left\{ \begin{array}{ll}
 \max\left[ \kappa_{tx1} v_k^{\max},
      v(x_k^+,s_k^+) + \kappa_{tx2}( v(x_k,s_k) - v(x_k^+,s_k^+))  \right]
                  & \mbox{if \eqref{eq:rhoc} hold,}\\
v_k^{\max} & \mbox{otherwise,}
\end{array}\right.
\end{eqnarray}
for some $\kappa_{tx1} \in (0,1)$ and $\kappa_{tx2}\in (0,1)$.

%%%%%%%%%%%%%%%%%%%%%%%%%%%%%%%%%%%%%%%%%%%%%%%%%%%%%%%%%%%%%%%%%%%%%%%%%%%%%%%%%%%%%%%%%

\subsubsection{Criticality step}\label{subsubsec:criticalitystep}

Two different criticality steps are employed: one for the subspaces $\calS_k$ with $dim(\calS_k)<n$ and one for the full space ($dim(\calS_k)=n$). In the latter, convergence is declared whenever at least one of the following conditions is satisfied: (1) the trust-region radius $\Delta_k$ is too small, (2) the computed direction $d_k$ is too small or (3) both feasibility and optimality have been achieved and the error between the real functions and the models is expected to be sufficiently small. As it was mentioned before, this error is directly linked to the $\Lambda$-poisedness measure given in Definition \ref{def:lambdapoisedness}. In the subspace, we are less demanding and only ask that either $\Delta_k$ be very small or both feasibility and optimality have been achieved without checking the models error though.

The complete criticality step in DEFT-FUNNEL is described in the next algorithm.

\algor{DFTRCRITSTEP}{CriticalityStep($\calS_k$, $\calY_k$, $\pi_{k-1}^f$, $\epsilon_i$)}{
%\vspace*{-0.3 cm}
\begin{description}
\item[Step 1:] If $dim(\calS_k)<0$, 
\begin{description}
\item[Step 1.1:] If $\Delta_k\leq \epsilon \|(x_k,s_k)\|$, return $(x_k,s_k)$.
\item[Step 1.2:] If $\|z(x_k,s_k)\| \leq \epsilon$ and $\hat{\pi}_i^f \leq \epsilon$, return $(x_k,s_k)$.
\end{description}
\end{description}

\begin{description}
\item[Step 2:] If $dim(\calS_k)=n$, 
\begin{description}
\item[Step 2.1:] If $\Delta_k\leq \epsilon \|(x_k,s_k)\|$ or $\|d_k\|\leq \epsilon \|(x_k,s_k)\|$, return $(x_k,s_k)$.
\item[Step 2.2:] Define $\hat{m}^f_i = m^f_k$, $\hat{m}^c_i = m^c_k$ and $\hat{\pi}_i^f=\pi_{k-1}^f$.
\item[Step 2.3:] If $\|z(x_k,s_k)\| \leq \epsilon_i$ and $\hat{\pi}_i^f \leq \epsilon_i$, set $\epsilon_{i+1}=\max\left[\alpha\|z(x_k,s_k)\|,\alpha\hat{\pi}_i^f,\epsilon\right]$ and modify $\calY_k$ as needed to ensure it is $\Lambda$-poised in $\calB(x_k;\epsilon_{i+1})$. If $\calY_k$ was modified, compute new models $\hat{m}_i^f$ and $\hat{m}_i^c$, calculate $\hat{r}_i$ and $\hat{\pi}_i^f$ and increment $i$ by one. If $\|z(x_k,s_k)\|\leq \epsilon$ and $\hat{\pi}_i^f\leq \epsilon$, return $(x_k,s_k)$; otherwise, start Step~2.3 again; 
\item[Step 2.4:] Set $m^f_k = \hat{m}^f_i$, $m^c_k = \hat{m}^c_i$, $\pi_{k-1}^f=\hat{\pi}_i^f$, $\Delta_k = \beta\max\left[\|z(x_k,s_k)\|,\pi_{k-1}^f\right]$ and define $\vartheta_i=x_k$ if a new model has been computed. 
\end{description}
\end{description}
}

%%%%%%%%%%%%%%%%%%%%%%%%%%%%%%%%%%%%%%%%%%%%%%%%%%%%%%%%%%%%%%%%%%%%%%%%%%%%%%%%%%%%%%%%%

\subsection{Maintenance of the interpolation set and trust-region updating strategy}\label{subsubsec:mainteninterset}

The management of the geometry of the interpolation set is based on the self-correcting geometry scheme proposed in \cite{ScheinbergToint10}, where unsuccessful trial points are used to improve the geometry of the interpolation set. It depends on the criterion used to define successful iterations, which is passed to the algorithm through the parameter \textit{criterion}. This parameter depends on the iteration type ($\mu$, $f$ or $z$, as explained previously). In general, unsuccessful trial points replace other sample points in the interpolation set which maximize a combined criteria of distance and poisedness involving the trial point. Finally, we also notice that we do not make use of ``dummy'' interpolation points resulting from projections onto the subspaces as in \cite{Sampaio2016} anymore. The whole procedure is described in Algorithm~\ref{DFTRINTSET}. The definition of the maximum cardinality of the interpolation set, $p^{\max}$, is given in Section~\ref{subsubsec:surr}.

%\vspace*{0.3 cm}
\algor{DFTRINTSET}{UpdateInterpolationSet($\calY_k$, $x_k$, $x_k^+$, 
$\Delta_k$, $\epsilon_i$, \textit{criterion})}{
%\vspace*{-0.3 cm}
\begin{description}
\item[1: Augment the interpolation set.] If $|\calY_k| < p^{\max}$, 
then define $\calY_{k+1}=\calY_k\cup\{x_k^+\}$.

\item[2: Successful iteration.] If $|\calY_k| = p^{\max}$ and \textit{criterion}
holds, then define $\calY_{k+1}=\calY_k\setminus\{y^{k,r}\}\cup\{x_k^+\}$ for 
\begin{eqnarray}
y^{k,r}=\arg\max_{y^{k,j}\in\calY_k}\|y^{k,j}-x_k^+\|^2|\ell_{k,j}(x_k^+)|.
\end{eqnarray}

\item[3: Replace a far interpolation point.] If $|\calY_k| = p^{\max}$, \textit{criterion} fails, either $x_k\neq \vartheta_i$ or $\Delta_k\leq\epsilon_i$, and the set 
\begin{eqnarray}
\calF_k
\eqdef\{y^{k,j}\in\calY_k\text{ such that }\|y^{k,j}-x_k\
|>\zeta\Delta \text{ and } \ell_{k,j}(x_k^+)\neq 0\}
\end{eqnarray}
is non-empty, then define $\calY_{k+1}=\calY_k\setminus\{y^{k,r}\}\cup\{x_k^+\}$, where 
\begin{eqnarray}
y^{k,r}=\arg\max_{y^{k,j}\in\calF_k}\|y^{k,j}-x_k^+\|^2|\ell_{k,j}(x_k^+)|.
\end{eqnarray}

\item[4: Replace a close interpolation point.] If $|\calY_k| = p^{\max}$, \textit{criterion} fails, either $x_k\neq \vartheta_i$ or $\Delta_k\leq\epsilon_i$, the set $\calF_k$ is empty, and the set 
\begin{eqnarray}
\calC_k\eqdef\{y^{k,j}\in\calY_k\text{ such that }\|y^{k,j}-x_k\|\leq\zeta\Delta \text{ and } |\ell_{k,j}(x_k^+)|>\Lambda\}
\end{eqnarray}
is non-empty, then define $\calY_{k+1}=\calY_k\setminus\{y^{k,r}\}\cup\{x_k^+\}$, where 
%$r$ is the index for any point in $\calC_k$, for instance, such that 
\begin{eqnarray}
y^{k,r}=\arg\max_{y^{k,j}\in\calC_k}\|y^{k,j}-x_k^+\|^2|\ell_{k,j}(x_k^+)|.
\end{eqnarray}

\item[5: No replacements.] If $|\calY_k| = p^{\max}$, \textit{criterion} fails and either [$x_k=\vartheta_i$ and $\Delta_k > \epsilon_i$] or $\calF_k\cup\calC_k=\emptyset$, then define $\calY_{k+1}=\calY_k$.
\end{description}
}

The trust-region update strategy associated to $f$- and $z$-iterations are now described. Following the idea proposed in~\cite{GrattonTointTroltzsch2011}, the trust-region radii are allowed to decrease even when the interpolation set has been changed after the replacement of a far point or a close point at unsuccessful iterations. However, the number of times it can be shrunk in this case is limited to $\nu_f^{\max}$ and $\nu_z^{\max}$ as a means to prevent the trust regions from becoming too small. If the interpolation set has not been updated, the algorithm understands that the lack of success is not due to the surrogate models but rather to the trust region size and thus it reduces the latter.

%\vspace*{0.3 cm}
\algor{DFTRFITER}{$f$-iteration($x_k$, $s_k$, $x_k^{+}$, $s_k^{+}$, $\Delta_k^f$, $\Delta_k^z$)}{
%\vspace*{-0.3 cm}
\begin{description}
\item[1: Successful iteration.] If $\rho_k^f\geq \eta_1$, set $(x_{k+1},s_{k+1}) = (x_k^{+},s_k^{+})$ and $\nu_f = 0$. The radius of $\calT_k$ is updated by
\begin{eqnarray}\label{eq:radupdatefsuc}
\Delta_{k+1}^f = \left\{ \begin{array}{ll}
  \min\left[ \max[ \gamma_2\|d_k\|,\Delta_k^f],\Delta^{\max}\right] & \tim{if} \rho_k^f \geq \eta_2, \\
  \Delta_k^f & \tim{if} \rho_k^f \in [\eta_1,\eta_2),\\
  \end{array} \right.
\end{eqnarray}
The radius of $\calN_k$ is updated by
\begin{eqnarray}\label{eq:radupdatecfsuc}
 \Delta_{k+1}^z = \left\{ \begin{array}{ll}
 \min\left[ \max[ \gamma_2\|n_k\|, \Delta_k^z],\Delta^{\max}\right] & 
 \!\!\!\tim{if} v(x_k^{+},s_k^{+}) < \eta_3\, v_k^{\max}, \\
 \Delta_k^z & \!\!\!\tim{otherwise.}
\end{array} \right.
\end{eqnarray}

\item[2: Unsuccessful iteration.] If $\rho_k^f<\eta_1$, set $(x_{k+1},s_{k+1}) = (x_k,s_k)$ and $\Delta^z_{k+1}=\delta^z_k$. The radius of $\calT_k$ is updated by
\begin{eqnarray}\label{eq:radupdatefunsuc} 
\Delta_{k+1}^f = \left\{ \begin{array}{ll}
\gamma_1 \|d_k\| & \!\!\!\!\!\tim{if either} (\calY_{k+1} \neq \calY_k \tim{and} \nu_f \leq \nu_f^{\max})\\
    & \,\;  \tim{or} \calY_{k+1} = \calY_k, \\
\Delta_k^f & \!\!\!\!\!\tim{if} \calY_{k+1} \neq \calY_k \tim{and} \nu_f > \nu_f^{\max}, \\
\end{array} \right. 
\end{eqnarray}
If $\calY_{k+1} \neq \calY_k$ and $\nu_f \leq \nu_f^{\max}$, update $\nu_f = \nu_f + 1$.
\end{description}
}

The operations related to $z$-iterations follow below. 

%\vspace*{0.3 cm}
\algor{DFTRCITER}{$z$-iteration($x_k$, $s_k$, $x_k^{+}$, $s_k^{+}$, $\Delta_k^f$, $\Delta_k^z$)}{
%\vspace*{-0.3 cm}
\begin{description}
\item[1: Successful iteration.] If \eqref{eq:rhoc} holds, set $(x_{k+1},s_{k+1}) = (x_k^{+},s_k^{+})$, $\Delta_{k+1}^f=\Delta_k^f$ and $\nu_z = 0$. The radius of $\calN_k$ is updated by

\begin{eqnarray}\label{eq:radupdatecsuc}
 \Delta_{k+1}^z = \left\{ \begin{array}{ll}
\min\left[ \max[ \gamma_2\|n_k\|,\Delta_k^z],\Delta^{\max}\right]
     		& \!\!\!\!\! \tim{if} \rho_k^z \geq \eta_2, \\
\Delta_k^z  & \!\!\!\!\! \tim{if} \rho_k^z \in [\eta_1,\eta_2). \\
  \end{array} \right. 
\end{eqnarray}

\item[2: Unsuccessful iteration.] If \eqref{eq:rhoc} fails, set $(x_{k+1},s_{k+1}) = (x_k,s_k)$ and $\Delta^f_{k+1}= \Delta^f_k$. The radius of $\calN_k$ is updated by

\begin{eqnarray}\label{eq:radupdatecunsuc} 
\Delta_{k+1}^z = \left\{ \begin{array}{ll}
\gamma_1 \|n_k\| & \!\!\!\!\!\tim{if} \|n_k\|\neq 0 \tim{and either} (\calY_{k+1} \neq \calY_k \tim{and} \nu_z \leq \nu_z^{\max})\\
    & \,\;  \tim{or} \calY_{k+1} = \calY_k, \\
\Delta_k^z & \!\!\!\!\!\tim{if} \calY_{k+1} \neq \calY_k \tim{and} \nu_z > \nu_z^{\max}, \\
\gamma_1 \Delta_k^z & \!\!\!\!\!\tim{if} \|n_k\|= 0, \\
\end{array} \right. 
\end{eqnarray}
If $\calY_{k+1} \neq \calY_k$ and $\nu_z \leq \nu_z^{\max}$, update $\nu_z = \nu_z + 1$.

\end{description}
}

%%%%%%%%%%%%%%%%%%%%%%%%%%%%%%%%%%%%%%%%%%%%%%%%%%%%%%%%%%%%%%%%%%%%%%%%%%%%%%%%%%%%%%%%%

\subsubsection{The local-search algorithm}

We now provide the full description of the local search which assembles all the previous subroutines.

%\vspace*{0.3 cm}
\algor{LocalSearch}{LocalSearch($\calS$, $\calY$, $x$, $s$, $\Delta^f$, $\Delta^z$, $v^{\max}$)}{
\begin{description}
\item[0: Initialization.]
Choose an initial vector of Lagrange multipliers
$\mu_{-1}$ and parameters $\epsilon>0$, $\epsilon_0>0$, $\Delta^f_0>0$, $\Delta^z_0>0$, $\alpha \in (0,1)$, $0 < \gamma_1 < 1 < \gamma_2$, $\zeta\geq 1$, $0 < \eta_1 < \eta_2 < 1$, $\Lambda>1$, $\beta>0$, $\eta_3 >0$. Define $\Delta_0=\min[\Delta^f_0,\Delta^z_0] \leq \Delta^{\max}$. Initialize $\calY_0$, with $x_0\in\calY_0 \subset \calB(x_0;\Delta_0)$ and $|\calY_0|\geq n+1$, as well as the maximum number of interpolation points $p_{\max}\geq |\calY_0|$. Compute the associated models $m_0^{f}$ and $m_0^{c}$ around $x_0$ and Lagrange polynomials $\{l_{0,j}\}_{j=0}^p$. Set $\calS_{-1}=\mathbb{R}^n$. Define $v_0^{\max} = \max[\kappa_{za}, \kappa_{zr} v(x_0,s_0)]$, where $\kappa_{za}>0$ and $\kappa_{zr}>1$. Compute $r_{-1}$ by solving \eqref{eq:optmeasureprob} with normal step $n_{-1}=0$ and define $\pi_{-1}^f$ as in \eqref{eq:optmeasuredef}. Define $\nu_f^{\max}>0$ and $\nu_z^{\max}>0$ and set $\nu_f=\nu_z=0$. Set $k=0$ and $i=0$.

\item[1:] \underline{SubspaceMinimization}($\calS_{k-1}$, $\calY_k$, $x_k$, $s_k$, 
$\Delta_k^f$, $\Delta_k^z$, $v_k^{\max}$).

\item[2:] \underline{CriticalityStep}($\calS_k$, $\calY_k$, $\pi_{k-1}^f$, $\epsilon_i$).

\item[3:] \underline{NormalStep}($x_k$, $s_k$, $\pi_k^v$, $v_k$, $v_k^{\max}$).

\item[4:] \underline{TangentStep}($x_k$, $s_k$, $n_k$).

\item[5: Conclude a $\mu$-iteration.]
If $n_k=t_k=0$, then
\vspace*{-2mm}
\begin{description}
\item[5.1:] set $(x_{k+1},s_{k+1})= (x_k,s_k)$, $\Delta_{k+1}^f = \Delta_k^f$ and
$\Delta_{k+1}^z = \Delta_k^z$;
\item[5.2:] set $\Delta_{k+1}=\min[\Delta^f_{k+1},\Delta^z_{k+1}]$, $v_{k+1}^{\max} = v_k^{\max}$ 
and $\calY_{k+1}=\calY_k$.
\end{description}
%\vspace*{-2mm}

\item[6: Conclude an $f$-iteration.]
If $t_k \neq 0$ and \eqref{eq:fmodratios} and \eqref{eq:feasok} hold,
\vspace*{-2mm}
\begin{description}
\item[6.1:] \underline{UpdateInterpolationSet}($\calY_k$ ,$x_k$, $x_k^{+}$, $\Delta_k$, $\epsilon_i$, `$\rho_k^f\geq \eta_1$');
\item[6.2:] \underline{$f$-iteration}($x_k$, $s_k$, $x_k^{+}$, $s_k^{+}$, $\Delta_k^f$, $\Delta_k^z$);
\item[6.3:] Set $\Delta_{k+1}=\min[\Delta^f_{k+1},\Delta^z_{k+1}]$ and $v_{k+1}^{\max} = v_k^{\max}$.
\end{description}

\item[7: Conclude a $z$-iteration.]
If either $n_k \neq 0$ and $t_k = 0$, or either one of \eqref{eq:fmodratios}
or \eqref{eq:feasok} fails, 

\vspace*{-2mm}
\begin{description}
\item[7.1:] \underline{UpdateInterpolationSet}($\calY_k$ ,$x_k$, $x_k^{+}$, $\Delta_k$, $\epsilon_i$, `\eqref{eq:rhoc}');
\item[7.2:] \underline{$z$-iteration}($x_k$, $s_k$, $x_k^{+}$, $s_k^{+}$, $\Delta_k^f$, $\Delta_k^z$);
\item[7.3:] Set $\Delta_{k+1}=\min[\Delta^f_{k+1},\Delta^z_{k+1}]$ and update $v_k^{\max}$ using \eqref{eq:ccmax}.
\end{description}

\item[8: Update the models and the Lagrange polynomials.] If $\calY_{k+1}\neq
  \calY_k$, compute the interpolation models $m_{k+1}^f$ and $m_{k+1}^c$
  around $x_{k+1}$ using $\calY_{k+1}$ and the associated Lagrange
  polynomials $\{l_{k+1,j}\}_{j=0}^p$. Increment $k$ by one and go to Step~1. 
\end{description}
}

\subsection{Parameters tuning and user goals}\label{subsec:paramstuning}

Like any other algorithm, DEFT-FUNNEL performance can be affected by how its parameters are tuned. In this section, we discuss about which parameters might have a major impact in its performance and how they should be tuned depending on the user goals. We consider four main aspects concerning the resolution of black-box problems: the budget, the priority level given to feasibility, the type of the objective function and the priority level given to global optimality. We can then describe the possible scenarios in decreasing order of difficulty as below:
\begin{itemize}
    \item \textbf{Budget:} very low, limited or unlimited.
    \item \textbf{Priority to feasibility:} high or low.
    \item \textbf{Objective function:} highly multimodal, multimodal or unimodal.
    \item \textbf{Priority given to global optimality:} high or low.
\end{itemize}

Clearly, if the objective function is highly multimodal and the budget is limited, one should give priority to a multistart strategy with a high number of samples per iteration in case where global optimality is important. This can be done via two ways: the first one is by reducing the budget for each local search through the optional input maxeval\_ls, which equals maxeval*0.7 by default; for instance, one could try setting maxeval\_ls = maxeval*0.4, so that each local search uses up to 40\% of the total budget only. The other possibility is to increase the number $N$ of random points generated at Step 3 in the global search (see Algorithm \ref{GlobalSearch}). In case where attaining the global minimum is not a condition, a better approach would be to give more budget to each local search so that the chances to reach a local minimum are higher. If the objective function is not highly multimodal, one should search for a good compromise between spending the budget on each local search and on the sampling of the multistart strategy.

We also notice that when the budget is too small and it is very hard to find a feasible solution, it may be a good idea to compute a tangent step only when feasibility has been achieved. This is due to the fact that tangent steps may still deteriorate the gains in feasibility obtained through normal steps. When this happens, more normal steps must be computed which requires more function evaluations. Therefore, instead of spending the budget with both tangent and normal steps without guarantee of feasibility in the end, it is a better strategy to compute only normal steps in the beginning so that the chances of obtaining a feasible solution are higher in the end. For this purpose, the user should set the constant value $\kappa_{\calR}=0$ in \eqref{eq:suffinside}. By doing so, the tangent step will be computed only if the normal step equals zero, which by \eqref{eq:normal_cond} it happens when the iterate is feasible. In the code, this can be done by setting \verb|kappa_r| to zero in the function \verb|deft_funnel_set_parameters|.

\section{Numerical experiments}\label{sec:numericalexp}

We divide the numerical experiments into two sections: the first one is focused on the evaluation of the performance of DEFT-FUNNEL on black-box optimization problems, while the second one aims at analyzing the benefits of the exploitation of white-box functions on grey-box problems. In all experiments with DEFT-FUNNEL, minimum $\ell_2$-norm models were employed. The criticality step threshold was set to $\epsilon=10^{-4}$ and the trust-region constants were defined as $\eta_1 = 10^{-4}$, $\eta_2 = 0.9$, $\eta_3 = 0.6$, $\gamma_1 =0.5$, $\gamma_2 = 2.0$, $\nu_f^{\max} = 20\times n$ and $\nu_z^{\max} = 20 \times n$, where $n$ is the number of variables.

We compare DEFT-FUNNEL with two popular algorithms for constrained black-box optimization: the Genetic Algorithm (GA) and the Pattern Search (PS) algorithm from the Matlab Global Optimization Toolbox \cite{Matlab2015}. In particular, GA has been set with the Adaptive Feasible ('mutationadaptfeasible') default mutation function and the Penalty ('penalty') nonlinear constraint algorithm. As for PS, since it is a local optimization algorithm, it has been coupled with the Latin Hypercube Sampling (LHS) \cite{McKay79} method in order to achieve global optimality, a common strategy found among practitioners. Our experiments cover therefore three of the most popular approaches for solving black-box problems: surrogate-based methods, genetic algorithms and pattern-search methods. Other DFO algorithms have already been compared against DEFT-FUNNEL in previous papers (see \cite{Sampaio2015,Sampaio2016}) on a much larger set of test problems mainly designed for local optimization.

\subsection{Black-box optimization problems}\label{subsec:bbprob}

The three methods are compared on a set of 14 well-known benchmark problems for constrained global optimization, including four industrial design problems --- Welded Beam Design (WB4) \cite{Deb2000}, Gas Transmission Compressor Design (GTCD4) \cite{Beightler1976}, Pressure Vessel Design (PVD4) \cite{Coello2002} and Speed Reducer Design (SR7) \cite{Floudas1990} --- and the Harley pooling problem (problem 5.2.2 from \cite{Floudas99}), which is originally a maximization problem and that has been converted into a minimization one. Besides the test problems originated from industrial applications, we have collected problems with different characteristics (multimodal, nonlinear, separable and non-separable, with connected and disconnected feasible regions) to have a broader view of the performance of the algorithms in various kinds of scenarios. For instance, the Hesse problem \cite{Hesse1973} is the result of the combination of 3 separable problems with 18 local minima and 1 global minimum, while the Gómez \#3 problem, listed as the third problem in \cite{Gomez1982}, consists of many disconnected feasible regions, thus having many local minima. The test problems G3-G11 are taken from the widely known benchmark list in \cite{Michalewicz1996}. In Table~\ref{table:bbprobs}, we give the number of decision variables, the number of constraints and the best known feasible objective function value of each test problem. 

\begin{table}[!ht]
\caption{Problem name, number of decision variables, number of constraints (simple bounds not included) and the best known feasible objective function value of each test problem.}
\centering
\begin{tabular}{ | c | c | c | c | } 
\hline
\multirow{2}{2cm}{\textbf{Test problem}} & \textbf{Number of} & \textbf{Number of} & \textbf{Best known feasible} \\ 
        & \textbf{variables} & \textbf{constraints} & \textbf{objective function value} \\ 
\hline
Harley (Harley Pooling Problem) & 9 & 6 & -600 \\ 
WB4 (Welded Beam Design) & 4 & 6 & 1.7250 \\ 
GTCD4 (Gas Transmission Compressor Design) & 4 & 1 & 2964893.85 \\ 
PVD4 (Pressure Vessel Design)  & 4 & 3 & 5804.45 \\ 
SR7 (Speed Reducer Design)  & 7 & 11 & 2994.42 \\ 
Hesse     & 6  & 6  & -310 \\ 
Gómez \#3 & 2  & 1  & -0.9711 \\ 
G3        & 2  & 1  & -1 \\ 
G4        & 5  & 6  & -30665.539 \\ 
G6        & 2  & 2  & -6961.8139 \\ 
G7        & 10 & 8  & 24.3062 \\ 
G8        & 2  & 2  & -0.0958 \\ 
G9        & 7  & 4  & 680.6301 \\ 
G11       & 2  & 1  & 0.75000455 \\ 
\hline
\end{tabular}
\label{table:bbprobs}
\end{table}

Two types of black-box experiments were conducted in order to compare the algorithms. In the first type, a small budget of 100 black-box calls is given to each algorithm in order to evaluate how they perform on difficult problems with highly expensive functions. In such scenarios, many algorithms have difficulties to obtain even local minima depending on the test problem. In the second type of experiments, we analyze their ability and speed to achieve global minima rather than local minima by allowing larger budgets that range from $100\times n$ to $500\times n$, where $n$ is the number of variables. We consider every function as black box in both types of experiments. Finally, we have run each algorithm 50 times on each test problem.

Only approximate feasible solutions of the problem \eqref{eq:originalproblem} are considered when comparing the best objective function values obtained by the algorithms, i.e. we require that each optimal solution $x^*$ satisfy $cv(x^*) \leq 10^{-4}$, where 
\begin{eqnarray}\label{eq:feassol}
cv \eqdef \max\left[ \left[z(x) - u^s\right]^+, \left[l^s - z(x)\right]^+ \right].
\end{eqnarray}

In the next two subsections, we show the results for the two types of experiments in the black-box setting.

\subsubsection{Budget-driven experiments}\label{subsubsec:firstexp}

The results of the first type of experiments are shown in Table \ref{table:test1}. In the second column, $f_{\textrm{OPT}}$ denotes the objective function value of the global minimum of the problem when it is known or the best known objective function value otherwise. For each solver, we show the best, the average and the worst objective function values obtained in 50 runs on every test problem.

\begin{table}[!ht]
\caption{Results for the first type of experiments on a budget of 100 black-box calls. For each solver, we show the best, the average and the worst objective function values obtained in 50 runs.}
\centering
\begin{tabular}{ccccccc}
\hline
\textbf{Prob} & \textbf{$f_{\textrm{OPT}}$} & \textbf{Solver} & \textbf{Best} & \textbf{Avg.} & \textbf{Worst} \\
\hline
Harley    & -600        & GA                   & -195.2264 & -3.1266  & 209.5655  \\
          &             & PS                   & None      & None     & None \\
          &             & \textbf{DEFT-FUNNEL} & -600      & -17.1508 & 301.9904 \\
\hline
WB4       & 1.7250      & GA                   & 2.3013 & 4.4525 & 6.1749  \\
          &             & PS                   & 3.6274 & 7.2408 & 11.0038 \\
          &             & \textbf{DEFT-FUNNEL} & 2.9246 & 7.3649 & 16.1476 \\
\hline 
GTCD4     & 2964893.85  & GA                   & 4004353.4045 & 9693274.3280  & 13860176.1776  \\
          &             & PS                   & 4953046.0849 & 12212940.1950 & 13786675.7772  \\
          &             & \textbf{DEFT-FUNNEL} & 3648033.5677 & 11080105.7502 & 28062428.3209  \\
\hline 
PVD4      & 5804.45     & GA                   & 5804.3762 & 5870.6299 & 6095.9716  \\
          &             & PS                   & 5877.6483 & 7650.8416 & 10827.5206 \\
          &             & \textbf{DEFT-FUNNEL} & 5804.3761 & 7360.2882 & 10033.0231 \\
\hline 
SR7       & 2994.42     & GA                   & 3027.8978 & 3151.6054 & 3438.7069  \\
          &             & PS                   & 3134.0525 & 3516.1761 & 5677.6224  \\
          &             & \textbf{DEFT-FUNNEL} & 3003.7577 & 3489.8743 & 4583.1364  \\
\hline 
Hesse     & -310        & GA                   & -292.1091 & -277.1073 & -259.1020  \\
          &             & PS                   & -302.1163 & -162.5650 & -49.4364   \\
          &             & \textbf{DEFT-FUNNEL} & -310      & -234.5969 & -24       \\
\hline 
Gómez \#3 & -0.9711     & GA                   & -0.9711 & -0.8551 & -0.6532  \\
          &             & PS                   & -0.9700 & -0.7774 & -0.4293  \\
          &             & \textbf{DEFT-FUNNEL} & -0.9711 & -0.0689 & 3.23333  \\
\hline 
G3        & -1          & GA                   & -1 & -0.9978 & -0.9872  \\
          &             & PS                   & -1 & -0.8162 & -0.0831  \\
          &             & \textbf{DEFT-FUNNEL} & -1 & -0.8967 & -0.0182  \\
\hline
G4        & -30665.539  & GA                   & -30674.0765 & -30048.4518 & -29243.6646  \\
          &             & PS                   & -30814.5885 & -29422.4521 & -27838.8541  \\
          &             & \textbf{DEFT-FUNNEL} & -31025.6056 & -30980.3524 & -29246.5608  \\
\hline
G6        & -6961.8139 & GA                    & -6240.7711 & -3520.1488 & -2275.7798  \\
          &             & PS                   & -6252.2652 & -3220.0072 & -1206.1356  \\
          &             & \textbf{DEFT-FUNNEL} & -6961.8165 & -6961.8158 & -6961.8146  \\
\hline
G7        & 24.3062091  & GA                   & 95.3300 & 325.3722 & 688.7635  \\
          &             & PS                   & 197.7475 & 434.6940 & 687.2492  \\
          &             & \textbf{DEFT-FUNNEL} & 24.3011 & 52.1011 & 185.5706  \\
\hline
G8        & -0.095825   & GA                   & -0.0909 & -0.0267 & -0.0002  \\
          &             & PS                   & -0.0953 & -0.0097 & 0.0185  \\
          &             & \textbf{DEFT-FUNNEL} & -0.0958 & -0.0471 & 0.0008  \\
\hline
G9        & 680.6300573 & GA                   & 775.1281 & 1253.1102 & 4815.6952  \\
          &             & PS                   & 953.2964 & 5336.9157 & 12000.4159 \\
          &             & \textbf{DEFT-FUNNEL} & 797.1996 & 1403.7992 & 2668.9271  \\
\hline
G11       & 0.75000455  & GA                   & 0.7500 & 0.7557 & 0.7941  \\
          &             & PS                   & 0.7502 & 0.8937 & 0.9997  \\
          &             & \textbf{DEFT-FUNNEL} & 0.7499 & 0.7513 & 0.8091  \\
\hline
\end{tabular}
\label{table:test1}
\end{table}

As it can be seen in Table \ref{table:test1}, DEFT-FUNNEL found the global minimum in 10 out of 14 problems, while GA and PS did it in only 5 problems. Besides, when considering the best value found by each solver, DEFT-FUNNEL was superior to the others or equal to the best solver in 12 problems. In the average and worse cases, it also presented a very good performance; in particular, its worse performance was inferior to all others' in only 4 problems. Although GA did not reach the global minimum often, it presented the best average-case performance among all methods, while PS presented the worst. Finally, PS was the only one unable to reach a feasible solution in the Harley pooling problem.

\subsubsection{Experiments driven by global minima}\label{subsubsec:secondexp}

We now evaluate the ability of each solver to find global minima rather than local minima. The results of the second type of experiments for each test problem are presented individually, which allows a better analysis of the evolution of each solver performance over the number of function evaluations allowed. Each figure is thus associated to one single test problem and shows the average progress curve of each solver after 50 trials as a function of the budget.

\begin{figure}
    \centering
    \begin{subfigure}[b]{0.475\textwidth}
        \centering
        \includegraphics[width=\textwidth]{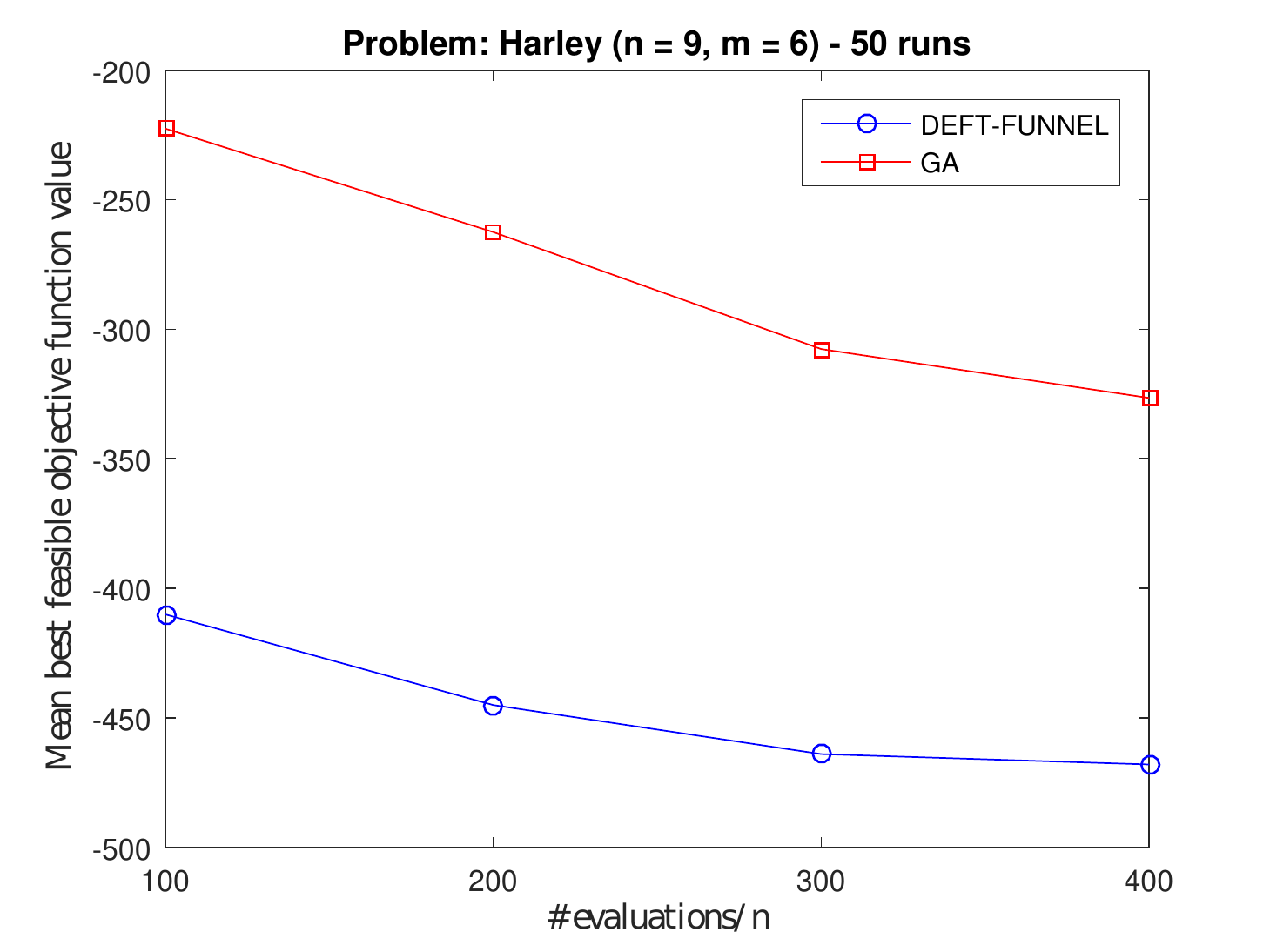}
        \caption{Test problem: Harley.}
        \label{fig:harley}
    \end{subfigure}
    \quad
    \begin{subfigure}[b]{0.475\textwidth}
        \centering
        \includegraphics[width=\textwidth]{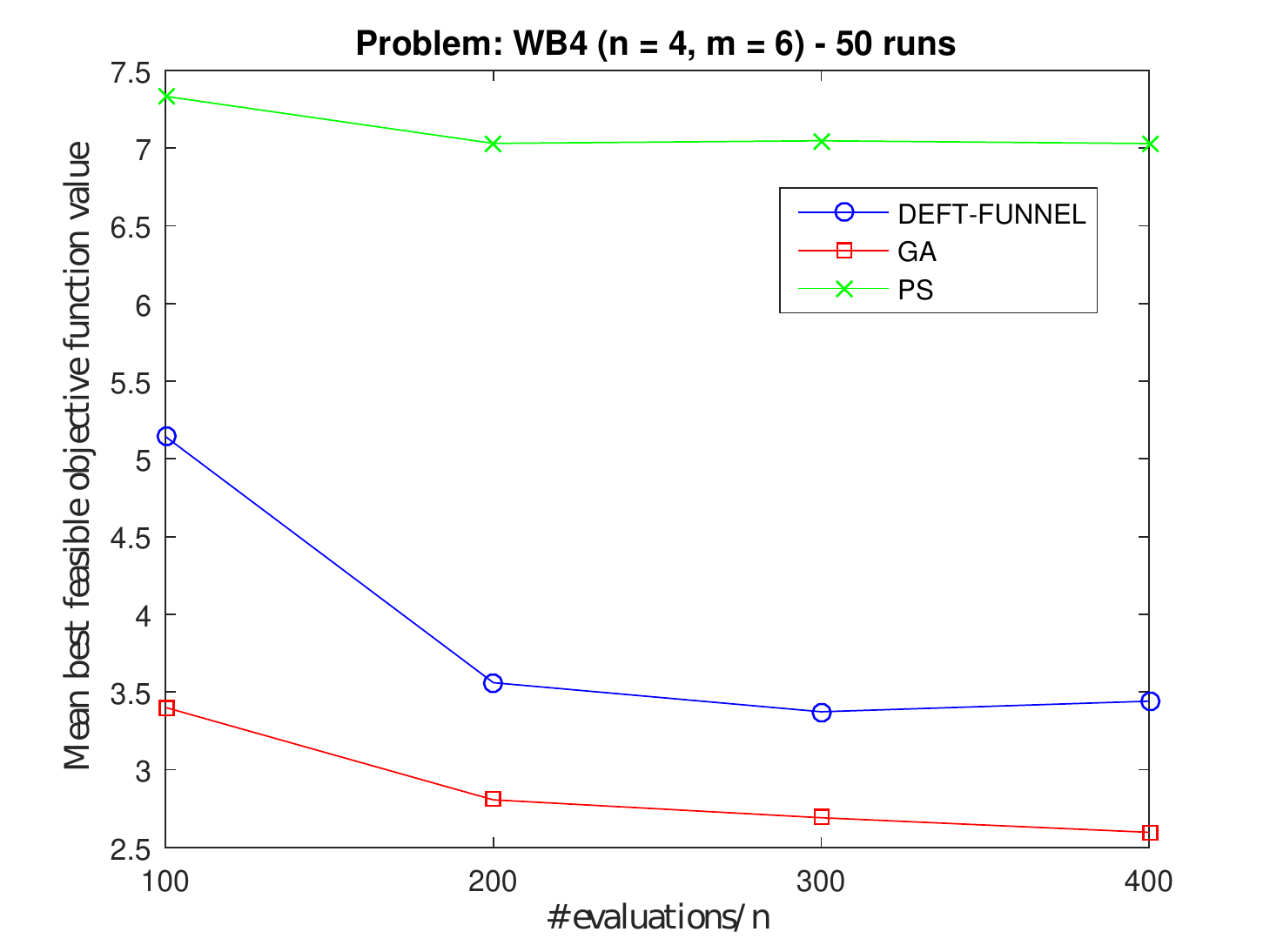}
        \caption{Test problem: WB4.}
        \label{fig:wb4}
    \end{subfigure}
    \vskip\baselineskip
    \begin{subfigure}[b]{0.475\textwidth}
        \centering
        \includegraphics[width=\textwidth]{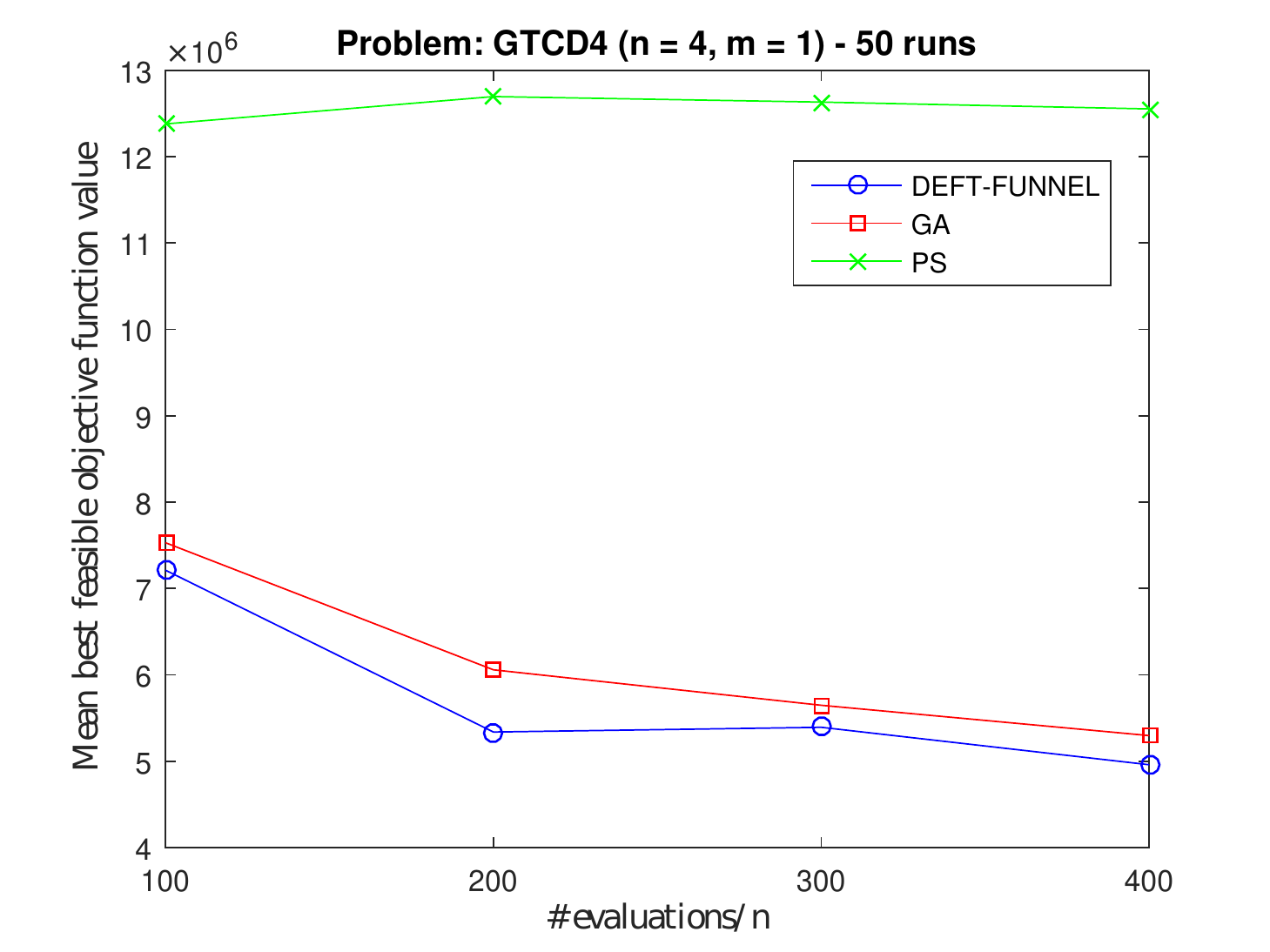}
         \caption{Test problem: GTCD4.}
        \label{fig:gtcd4}
    \end{subfigure}
    \quad
    \begin{subfigure}[b]{0.475\textwidth}
       \centering
        \includegraphics[width=\textwidth]{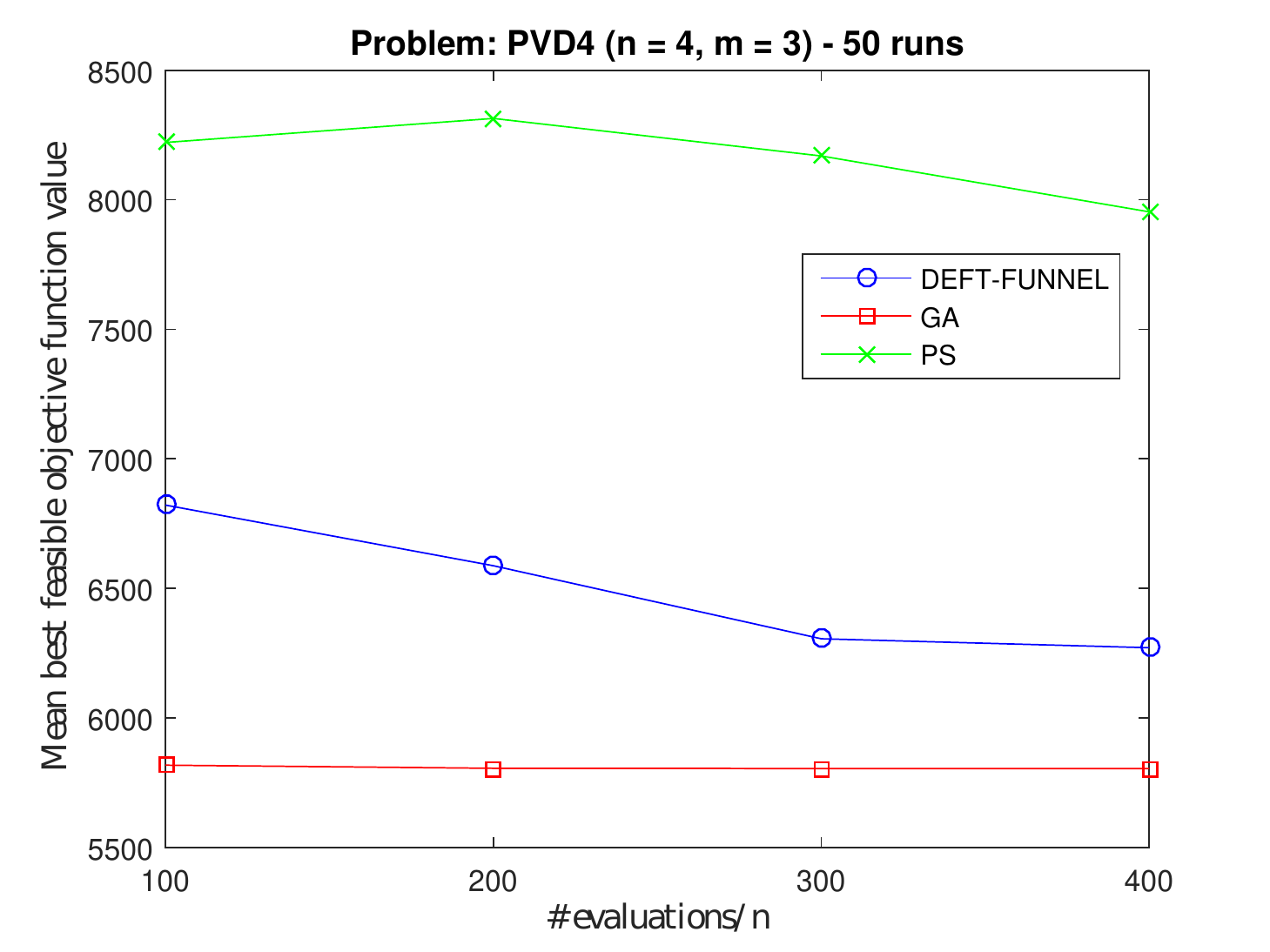}
        \caption{Test problem: PVD4.}
        \label{fig:pvd4}
    \end{subfigure}
    \caption{Mean best feasible objective function value obtained by each solver on 50 trials with budgets of $100\times n$, $200\times n$, $300\times n$ and $400\times n$ black-box function evaluations.}\label{fig:bbtests_globalopt_1}
\end{figure}

In Figure \ref{fig:bbtests_globalopt_1}, we show the results on test problems Harley, WB4, GTCD4 and PVD4. DEFT-FUNNEL and GA presented the best performance among the three methods, each one being superior to the other in 2 out of 4 problems. PS not only was inferior to the other two methods, but it also seemed not to be affected by the number of black-box calls allowed, as it can be seen in Figures \ref{fig:harley}, \ref{fig:wb4} and \ref{fig:gtcd4}. In particular, it did not find a feasible solution for the Harley problem at all. Moreover, the solutions found by PS had often much larger objective function values than those obtained by DEFT-FUNNEL and GA.

\begin{figure}
    \centering
    \begin{subfigure}[b]{0.475\textwidth}
       \centering
        \includegraphics[width=\textwidth]{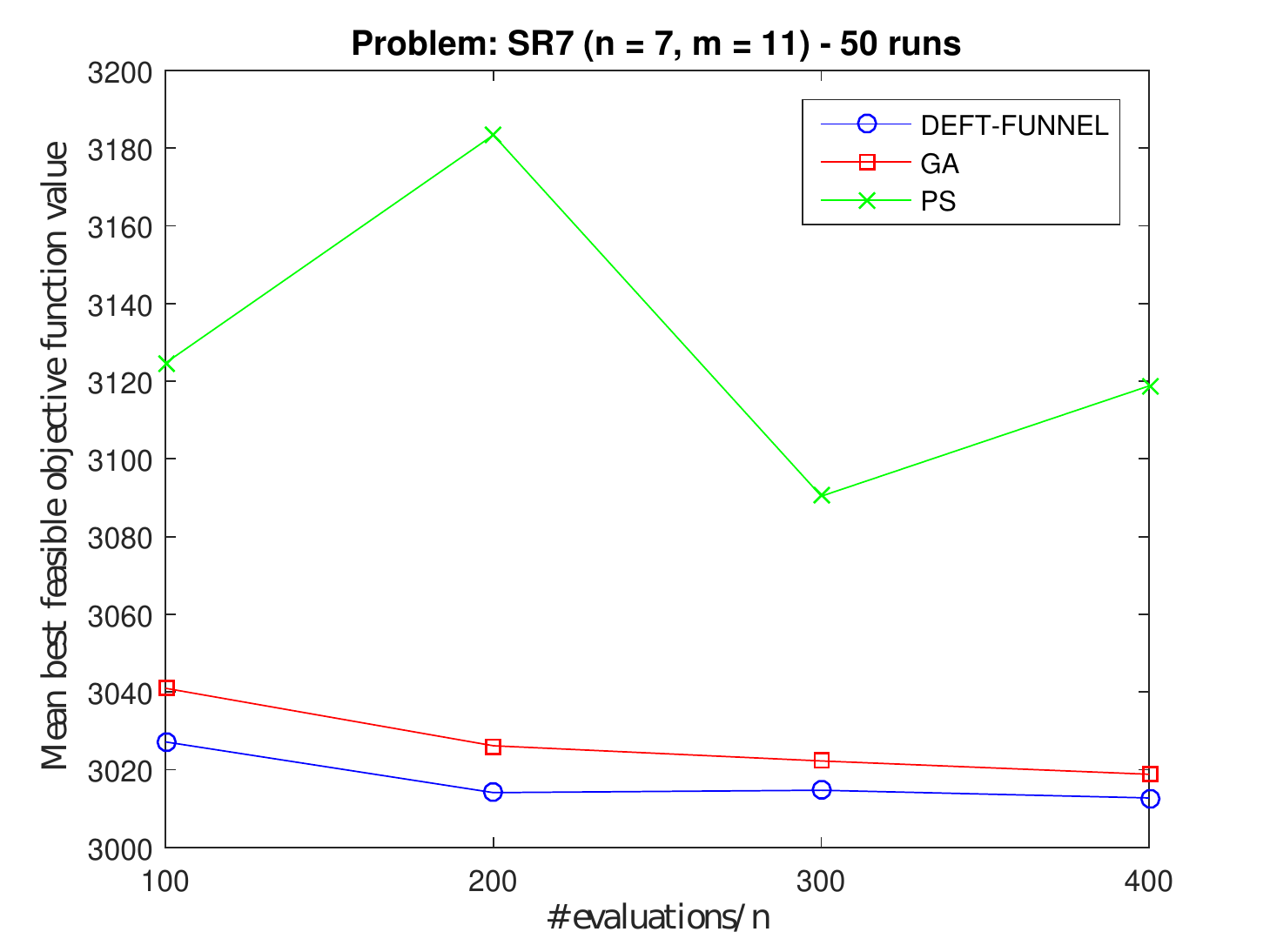}
        \caption{Test problem: SR7.}
        \label{fig:sr7}
    \end{subfigure}
    \quad
    \begin{subfigure}[b]{0.475\textwidth}
        \centering
        \includegraphics[width=\textwidth]{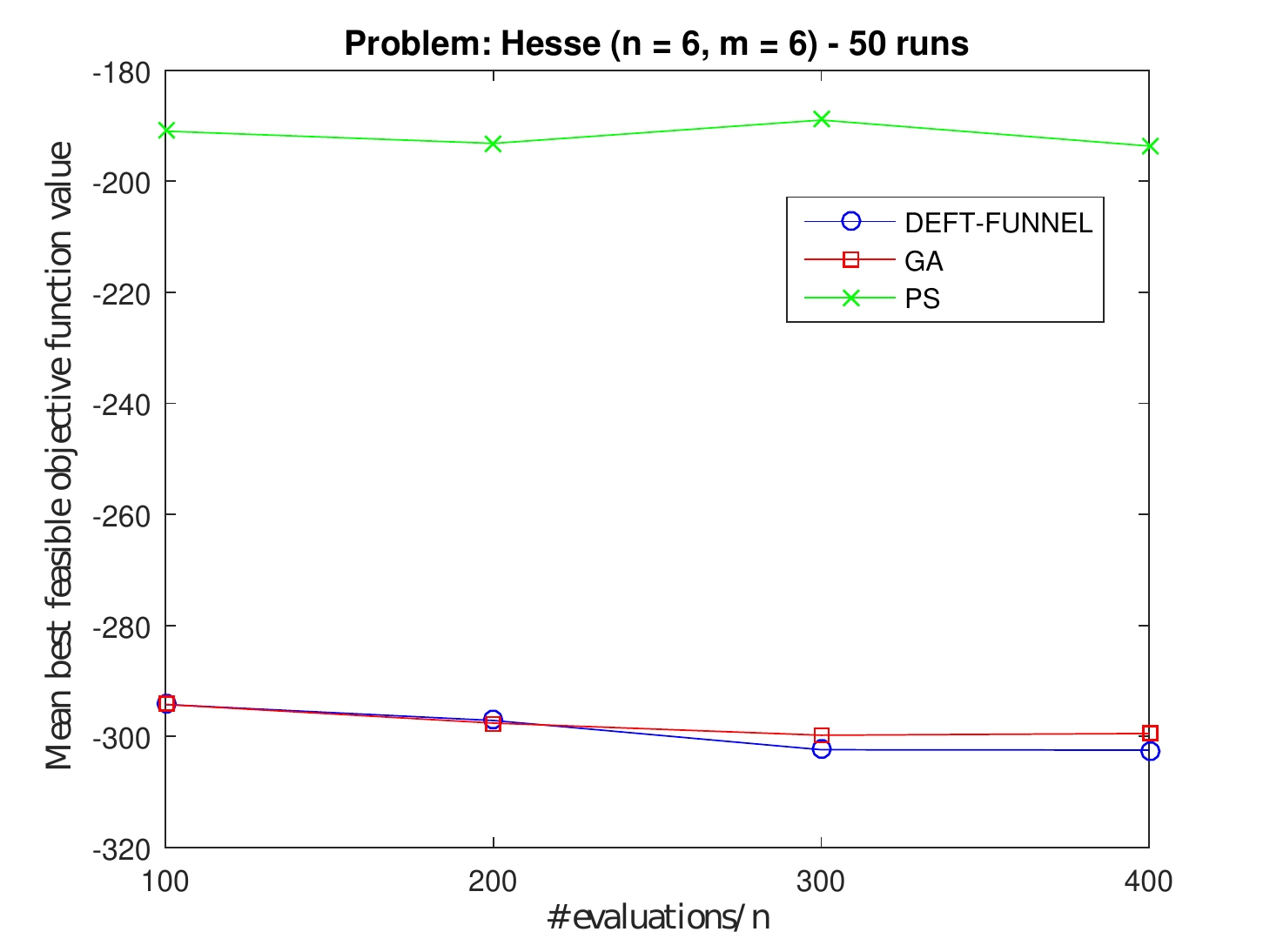}
        \caption{Test problem: Hesse.}
        \label{fig:hesse}
    \end{subfigure}
    \vskip\baselineskip
    \begin{subfigure}[b]{0.475\textwidth}
        \centering
        \includegraphics[width=\textwidth]{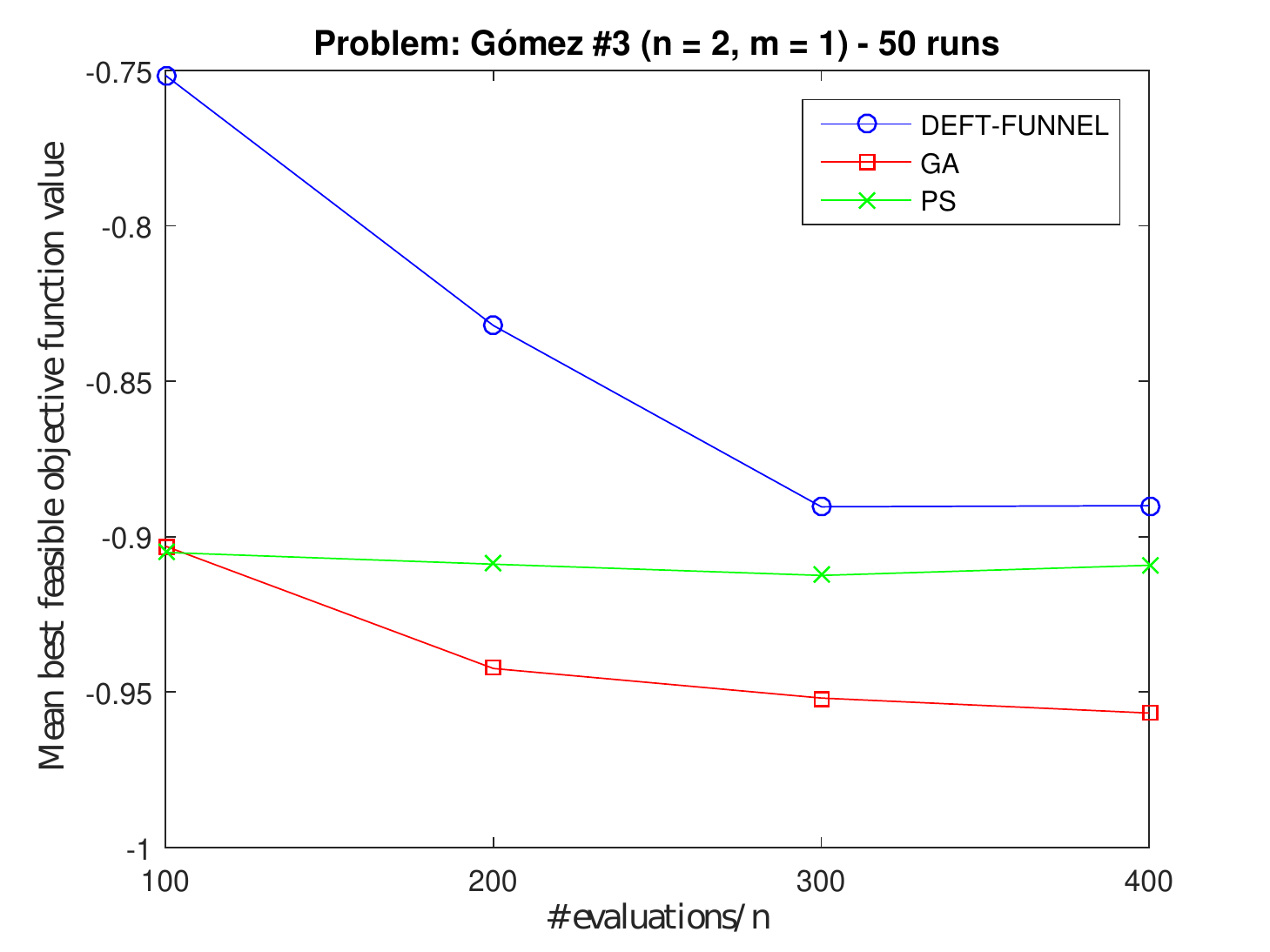}
         \caption{Test problem: Gómez \#3.}
        \label{fig:gomez3}
    \end{subfigure}
    \quad
    \begin{subfigure}[b]{0.475\textwidth}
       \centering
        \includegraphics[width=\textwidth]{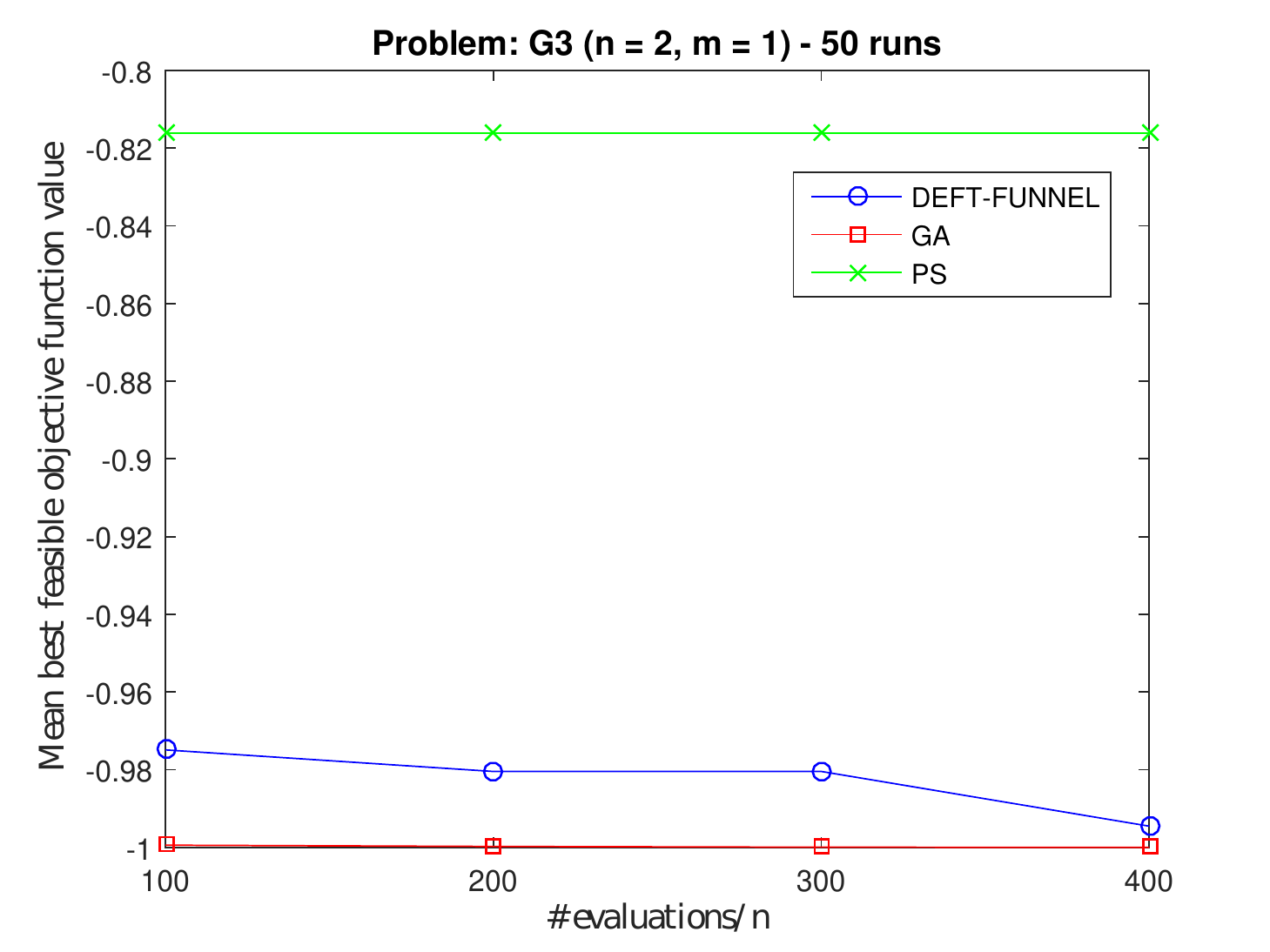}
        \caption{Test problem: G3.}
        \label{fig:g3}
    \end{subfigure}
    \caption{Mean best feasible objective function value obtained by each solver on 50 trials with budgets of $100\times n$, $200\times n$, $300\times n$ and $400\times n$ black-box function evaluations.}\label{fig:bbtests_globalopt_2}
\end{figure}

Figure \ref{fig:bbtests_globalopt_2} shows the results on test problems SR7, Hesse, Gómez \#3 and G3. DEFT-FUNNEL and GA had similar performances on SR7, Hesse and G3 problems, while PS had the poorest performance. GA was the best method on Gómez \#3, being very close to the global minimum in the end. Finally, PS performance on SR7 shows that allowing more black-box calls was not helpful and that the objective function value even increased.

\begin{figure}
    \centering
    \begin{subfigure}[b]{0.475\textwidth}
       \centering
        \includegraphics[width=\textwidth]{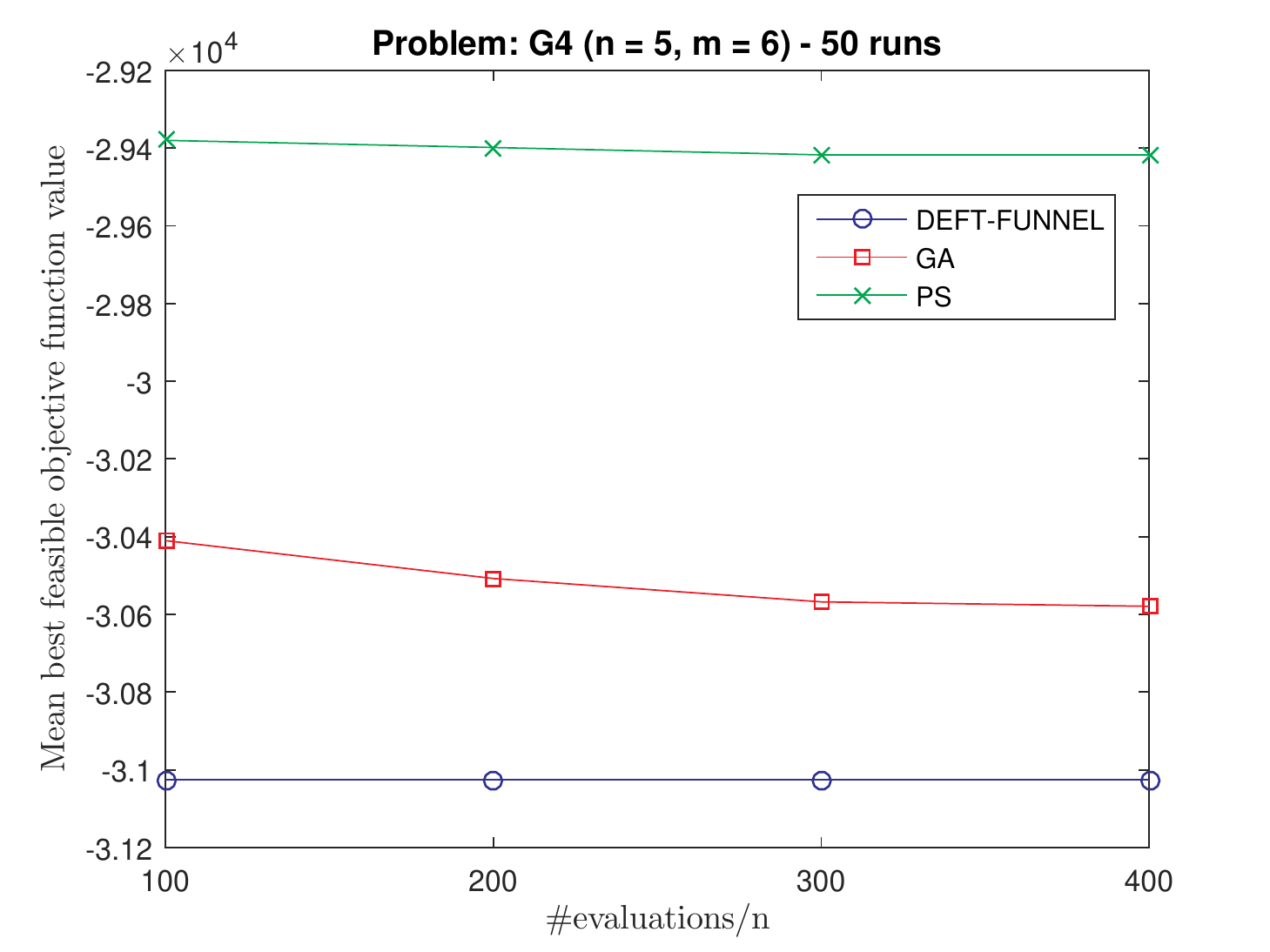}
        \caption{Test problem: G4.}
        \label{fig:g4}
    \end{subfigure}
    \quad
    \begin{subfigure}[b]{0.475\textwidth}
        \centering
        \includegraphics[width=\textwidth]{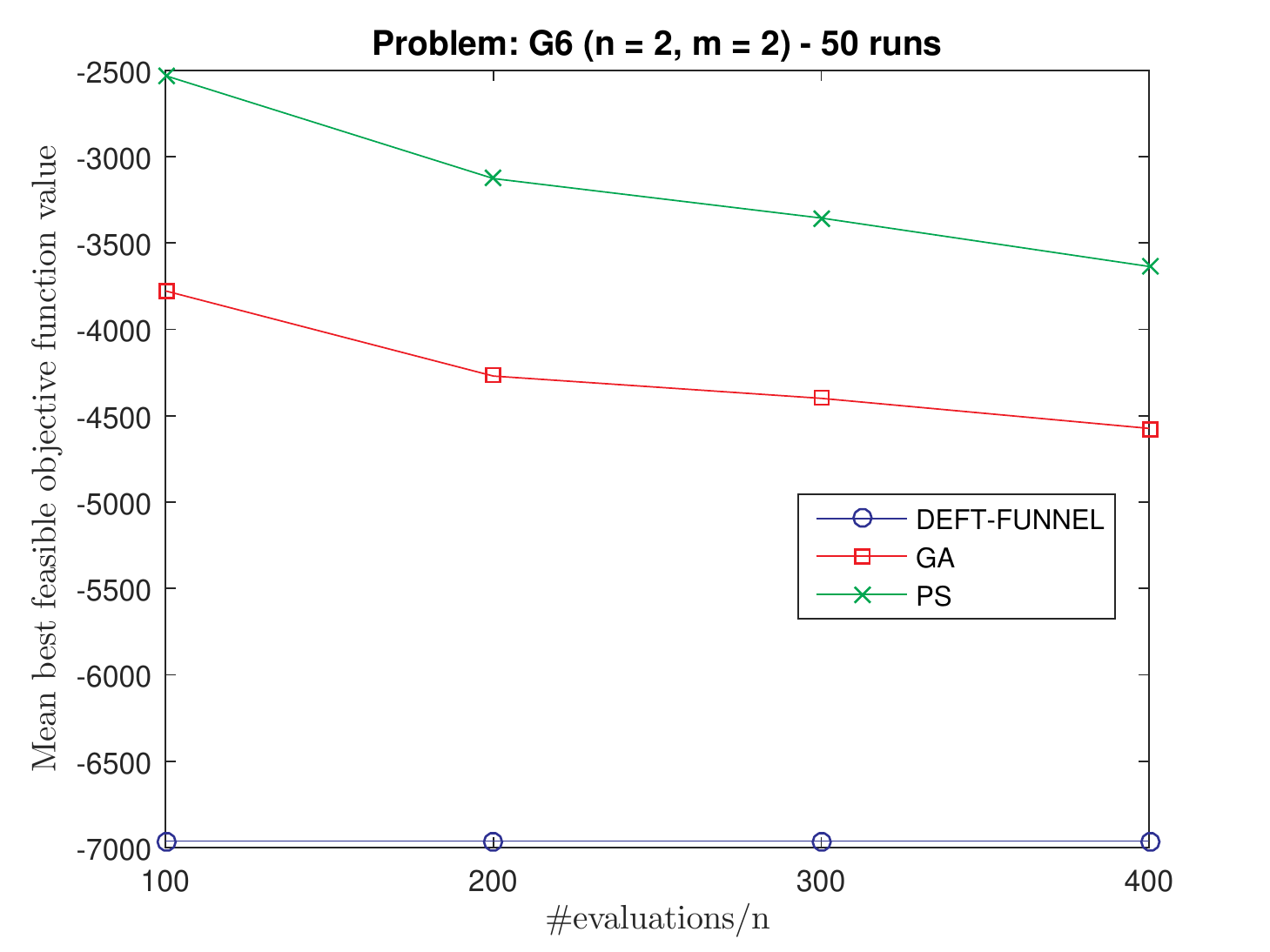}
        \caption{Test problem: G6.}
        \label{fig:g6}
    \end{subfigure}
    \vskip\baselineskip
    \begin{subfigure}[b]{0.475\textwidth}
        \centering
        \includegraphics[width=\textwidth]{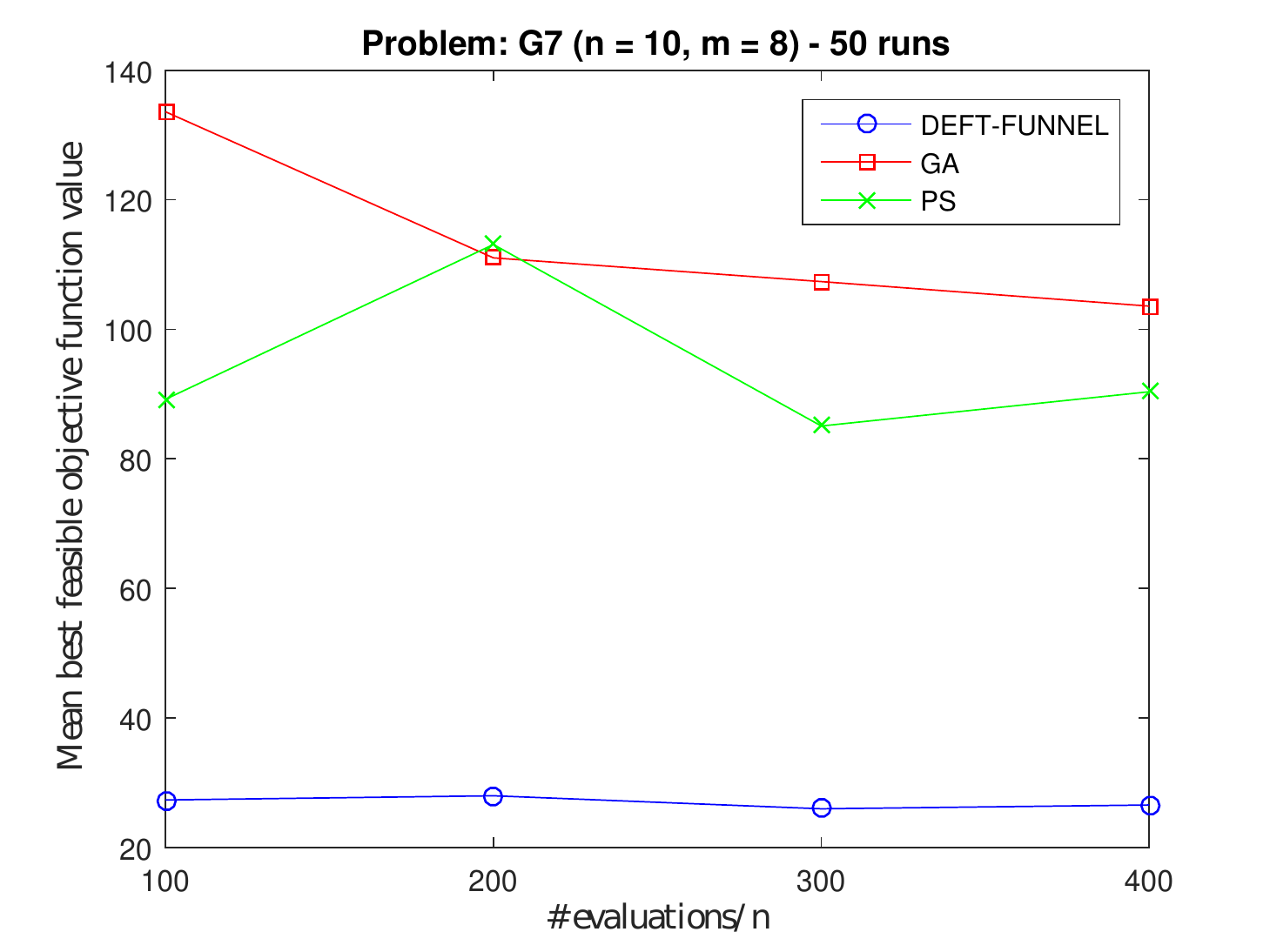}
         \caption{Test problem: G7.}
        \label{fig:g7}
    \end{subfigure}
    \quad
    \begin{subfigure}[b]{0.475\textwidth}
       \centering
        \includegraphics[width=\textwidth]{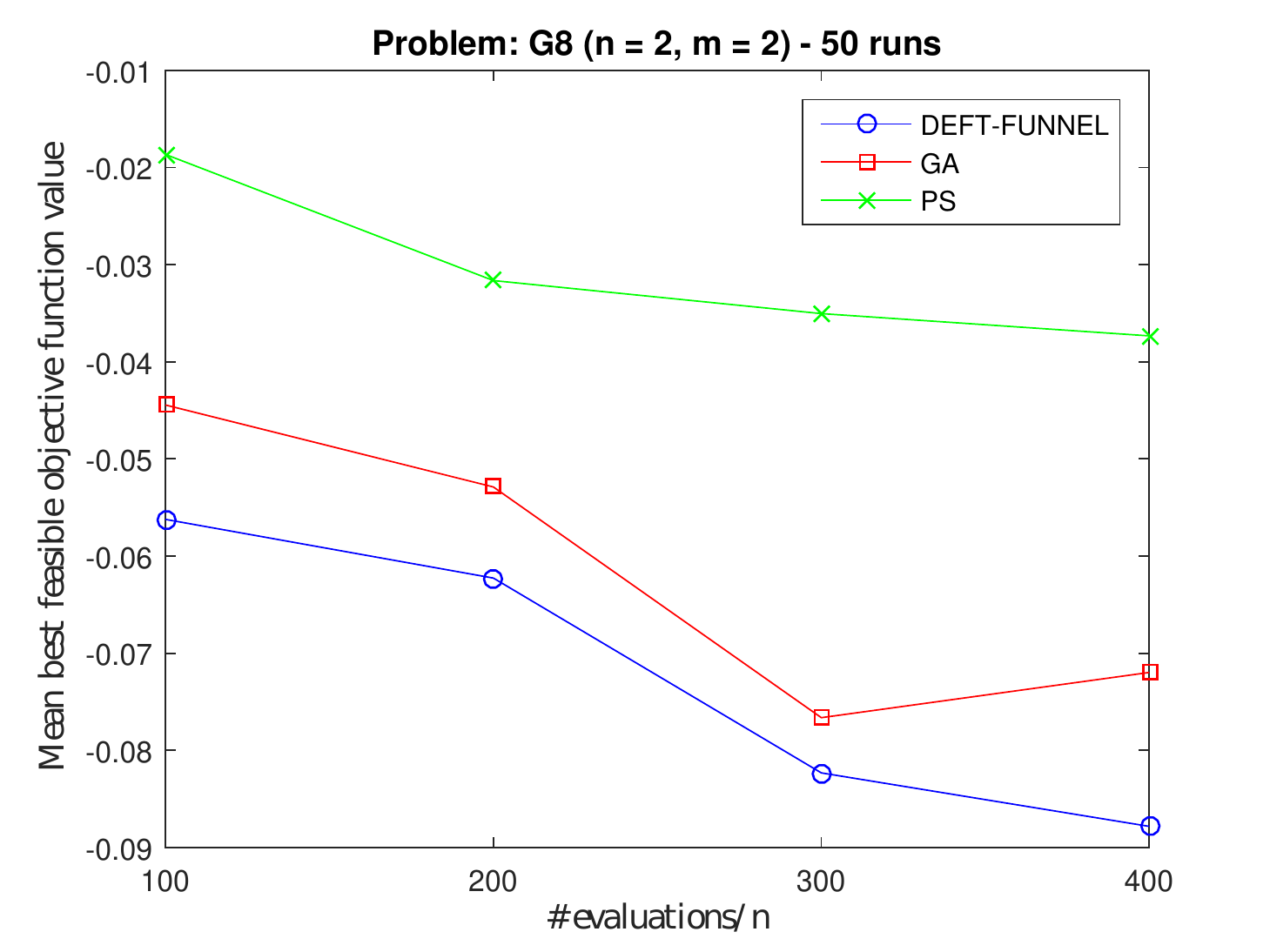}
        \caption{Test problem: G8.}
        \label{fig:g8}
    \end{subfigure}
    \vskip\baselineskip
    \begin{subfigure}[b]{0.475\textwidth}
       \centering
        \includegraphics[width=\textwidth]{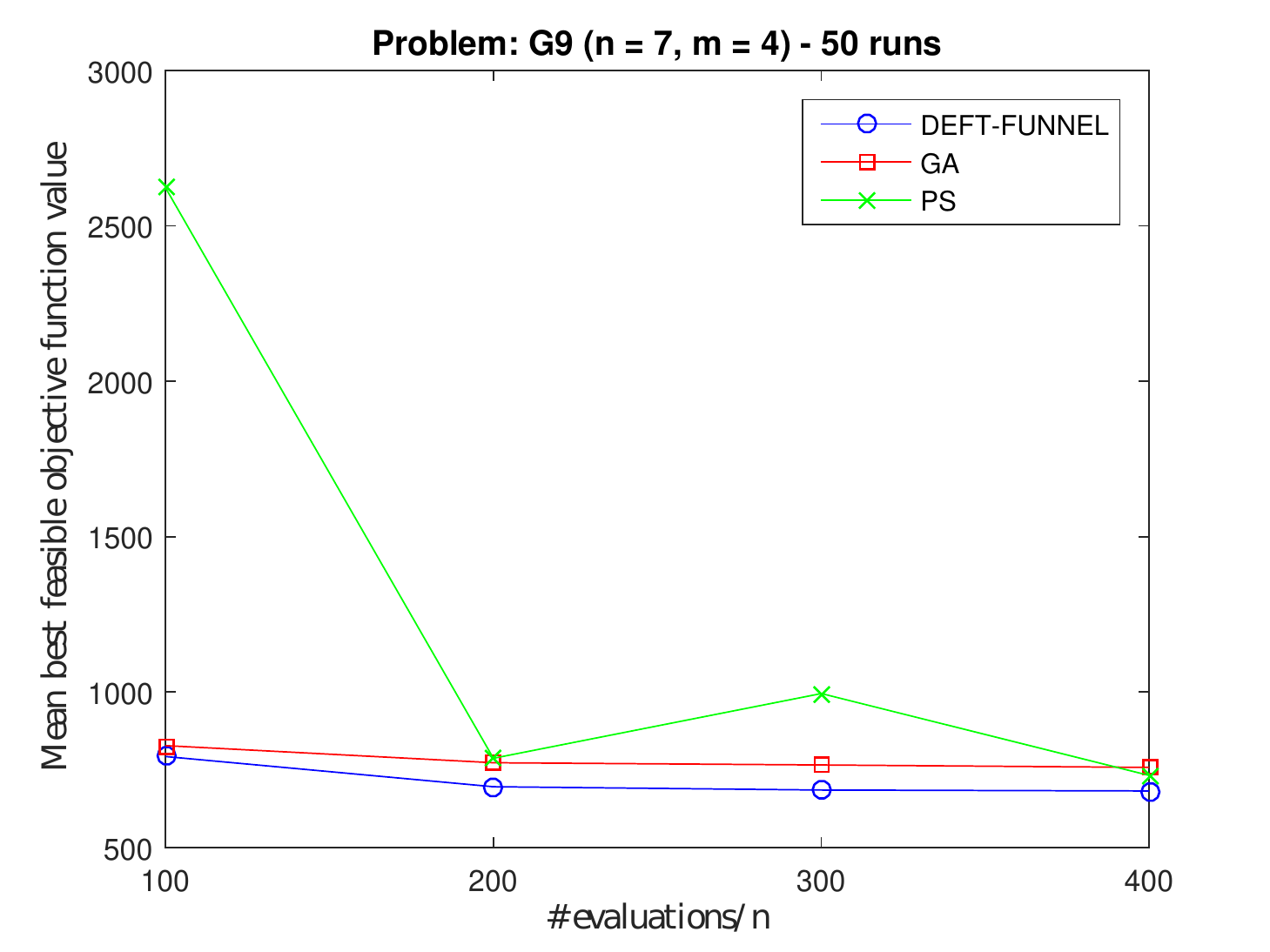}
        \caption{Test problem: G9.}
        \label{fig:g9}
    \end{subfigure}
    \quad
    \begin{subfigure}[b]{0.475\textwidth}
       \centering
        \includegraphics[width=\textwidth]{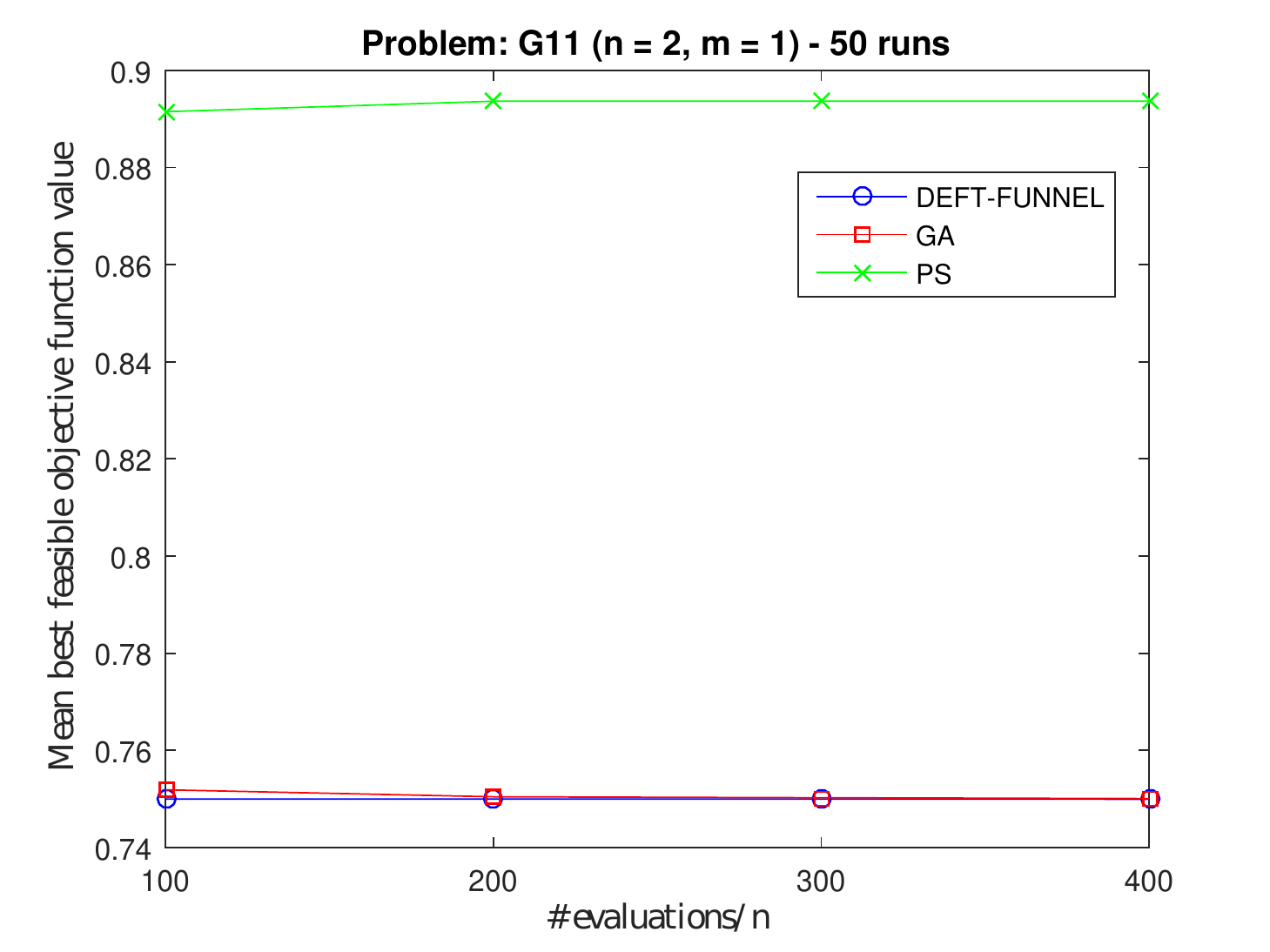}
        \caption{Test problem: G11.}
        \label{fig:g11}
    \end{subfigure}
    \caption{Mean best feasible objective function value obtained by each solver on 50 trials with budgets of $100\times n$, $200\times n$, $300\times n$ and $400\times n$ black-box function evaluations.}\label{fig:bbtests_globalopt_3}
\end{figure}

Figure \ref{fig:bbtests_globalopt_3} shows the results on test problems G4, G6, G7, G8, G9 and G11. DEFT-FUNNEL was superior to all other methods, attaining global optimality on 4 problems independently of the number of evaluations allowed. GA was the second best, followed by PS. Besides that PS did not find a feasible solution on Harley pooling problem, there is a significant increase on its objective function value as the number of evaluations grows on test problems SR7 and G7. These 3 problems are the ones with the largest number of constraints, which could be a reason for its poor performance.

\subsection{Grey-box optimization problems}\label{subsec:gbprob}

In this section, DEFT-FUNNEL is compared with GA and PS on a set of artificial grey-box problems built by selecting a test problem and defining some of its functions as black boxes and the remaining as white boxes. Among the 5 grey-box problems considered in the experiments, 3 were used in the black-box experiments described in the previous subsection, where all functions were considered as black boxes, namely: Hesse, SR7 and GTCD4. The reuse of these test problems allows for evaluating the performance improvement of DEFT-FUNNEL in cases where some information is available and for comparing its performance with that obtained in the black-box setting. The remaining grey-box problems are the problems 21 and 23 from the well-known Hock-Schittkowski collection \cite{Hock1980}, denoted here as HS21 and HS23, respectively. Both HS21 and HS23 have nonlinear objective functions; however, in HS21 the objective function is defined as white box, while in HS23 it is defined is black box. The only constraint present in HS21 is defined as black box, while in HS23 there is a balance between both categories among the constraints. We have therefore attempted to cover different possibilities related to the type of functions in a grey-box setting. The derivatives of all functions defined as white boxes were given as input to DEFT-FUNNEL. The reader can find more details about the tested grey-box problems in Table \ref{table:gbprobs}.

\begin{table}[!ht]
\caption{Problem name, number of decision variables, number of black-box constraints, number of white-box constraints, type of objective function -- black box (BB) or white box (WB) -- and the best known feasible objective function value of each test problem.}
\centering
\begin{tabular}{| c | c | c | c | c | c |} 
\hline
\multirow{2}{2cm}{\textbf{Test problem}} & \textbf{Number of} & \textbf{Number of} & \textbf{Number of} & \textbf{Type of} & \textbf{Best known feasible} \\ 
        & \textbf{variables} & \textbf{BB constraints} & \textbf{WB constraints} & \textbf{Obj. function} & \textbf{obj. function value} \\ 
\hline
GTCD4 & 4 & 0 & 1 & BB & 2964893.85 \\ 
SR7   & 7 & 9 & 2 & WB & 2994.42 \\ 
Hesse & 6 & 3 & 3 & WB & -310 \\ 
HS21  & 2 & 1 & 0 & WB & -99.96 \\ 
HS23  & 2 & 3 & 2 & WB  & 2 \\ 

\hline
\end{tabular}
\label{table:gbprobs}
\end{table}

In Table \ref{table:gbtests}, the results obtained with a budget of 100 black-box calls are presented. It can be seen that DEFT-FUNNEL was the only one to reach the global minimum in every problem, while GA had the best average-case and worst-case performances in general. We also notice that the three solvers had similar performances on HS21, attaining all the global minimum without much difficulty. 

When comparing the black-box and grey-box results obtained by DEFT-FUNNEL on problems GTCD4 and SR7, it is evident that the available information from the white-box functions contributed to improve its performance. Not only the best objective function values on these problems were improved, allowing for reaching the global minimum on SR7, but also the average and worst ones. Since DEFT-FUNNEL had already attained the global minimum on Hesse in the black-box setting, the only expected improvement would be in the average and worst cases, which did not happen in our experiments. This is probably due to the fact that this a multimodal problem with 18 local minima, being a combination of 3 separable problems and having 1 global minimum. Therefore, even if information about the function is partially available, the problem remains very difficult to solve due to its nature. 

\begin{table}[!ht]
\caption{Experiments with grey-box problems with a budget of 100 black-box calls. For each solver, we show the best, the average and the worst objective function values obtained in 50 runs.}
\centering
\begin{tabular}{ccccccc}
\hline
\textbf{Prob} & \textbf{$f_{\textrm{OPT}}$} & \textbf{Solver} & \textbf{Best} & \textbf{Avg.} & \textbf{Worst} \\
\hline
GTCD4     & 2964893.85  & GA                   & 4004353.4045 & 9693274.3280  & 13860176.1776  \\
          &             & PS                   & 4953046.0849 & 12212940.1950 & 13786675.7772  \\
          &             & \textbf{DEFT-FUNNEL} & 3564559.7818 & 9994027.4513  & 13937534.7346  \\
\hline
SR7       & 2994.42     & GA                   & 3027.8978 & 3151.6054 & 3438.7069  \\
          &             & PS                   & 3134.0525 & 3516.1761 & 5677.6224  \\
          &             & \textbf{DEFT-FUNNEL} & 2994.4244 & 3370.8758 & 4274.8433  \\
\hline
Hesse     & -310        & GA                   & -292.1091 & -277.1073 & -259.1020  \\
          &             & PS                   & -302.1163 & -162.5650 & -49.4364   \\
          &             & \textbf{DEFT-FUNNEL} & -310      & -210.6784 & -36   \\
\hline
HS21      & -99.96      & GA                   & -99.9599 & -99.5528 & -95.6688  \\
          &             & PS                   & -99.9599 & -99.9435 & -99.8086   \\
          &             & \textbf{DEFT-FUNNEL} & -99.96   & -99.9593 & -99.9267  \\
\hline
HS23      & 2           & GA                   & 4.7868 & 13.7529  & 174.1049  \\
          &             & PS                   & 4.5138 & 45.6059  & 274.6368  \\
          &             & \textbf{DEFT-FUNNEL} & 2      & 424.9427 & 1242.011  \\
\end{tabular}
\label{table:gbtests}
\end{table}

\section{Conclusions}\label{sec:conclusions}

We have introduced a new global optimization solver for general constrained grey-box and black-box optimization problems named DEFT-FUNNEL. It combines a stochastic multistart strategy with a trust-funnel sequential quadratic optimization algorithm where the black-box functions are replaced by surrogate models built from polynomial interpolation. The proposed solver is able to exploit available information from white-box functions rather than considering them as black boxes. Its code is open source and freely available at the Github repository \href{http://github.com/phrsampaio/deft-funnel}{http://github.com/phrsampaio/deft-funnel}. Unlike many black-box optimization algorithms, it can handle both inequality and equality constraints and it does not require feasible starting points. Moreover, the constraints are treated individually rather than grouped into a penalty function.

We have shown that DEFT-FUNNEL compares favourably with other state-of-the-art algorithms available for the optimization community on a collection of well-known benchmark problems. The reported numerical results indicate that the proposed approach is very promising for grey-box and black-box optimization. Future research works include extensions for multiobjective optimization and mixed-integer nonlinear optimization as well as as parallel version of the solver.

\bibliographystyle{abbrvnat}
\bibliography{references}

\end{document}